\documentclass[a4paper,12pt,reqno]{article}
\usepackage{pstricks, pst-plot, pst-node, pst-tree}
\usepackage{amsmath,amsfonts,amsthm,amssymb}

\def\EE{{\cal E}}
\def\HH{{\cal H}}
\newcommand{\bx}{{\bf x}}
\newcommand{\der}{\delta}
\newcommand{\dom}{\mbox{Dom}}
\newcommand{\ha}{\hat a}

\newcommand{\hr}{\hat r}
\newcommand{\hrk}{\hat r^{2\ka}}
\newcommand{\hz}{\hat z}
\newcommand{\ka}{\kappa}

\newcommand{\id}{\mbox{Id}}
\newcommand{\iou}{\int_{0}^{1}}
\newcommand{\iot}{\int_{0}^{t}}
\newcommand{\ist}{\int_{s}^{t}}
\newcommand{\norm}[1]{\lVert #1\rVert}

\newcommand{\ott}{[0,T]}

\newcommand{\xd}{{\bf x^{2}}}
\newcommand{\xdst}{{\bf x}_{st}^{\bf 2}}
\newcommand{\xds}{{\bf x^{2, s}}}
\newcommand{\xsst}{{\bf x}_{st}^{\bf 2,s}}
\newcommand{\1}{{\bf 1}}
\newcommand{\E}{{\rm E}}
\def\lpa{\langle}
\def\rpa{\rangle}
\def\sk{{\mathbb{D}}}

\newcommand{\R}{\mathbb R}
\newcommand{\N}{\mathbb N}

\newcommand{\cb}{\mathcal B}
\newcommand{\cac}{\mathcal C}

\newcommand{\cj}{\mathcal J}
\newcommand{\cl}{\mathcal L}
\newcommand{\cn}{\mathcal N}
\newcommand{\cq}{\mathcal Q}

\newcommand{\cz}{\mathcal Z}

\newcommand{\al}{\alpha}

\newcommand{\ga}{\gamma}
\newcommand{\la}{\lambda}
\newcommand{\laa}{\Lambda}

\newcommand{\si}{\sigma}

\newcommand{\vp}{\varphi}
\newcommand{\ze}{\zeta}

\newcommand{\lp}{\left(}
\newcommand{\rp}{\right)}
\newcommand{\lc}{\left[}
\newcommand{\rc}{\right]}
\newcommand{\lcl}{\left\{}
\newcommand{\rcl}{\right\}}
\newcommand{\lln}{\left|}
\newcommand{\rrn}{\right|}

\newcommand{\fin}
{ \vspace{-0.6cm}
\begin{flushright}
\mbox{$\Box$}
\end{flushright}
\noindent }

\newtheorem{theorem}{Theorem}[section]

\newtheorem{corollary}[theorem]{Corollary}

\newtheorem{definition}[theorem]{Definition}

\newtheorem{hypothesis}[theorem]{Hypothesis}
\newtheorem{lemma}[theorem]{Lemma}

\newtheorem{proposition}[theorem]{Proposition}

\newtheorem{remark}[theorem]{Remark}

\begin{document}

\title{Trees and asymptotic developments for fractional stochastic differential equations}
\author{
A. Neuenkirch\footnote{
Johann Wolfgang Goethe-Universit{\"a}t Frankfurt am Main,
FB Informatik und Mathematik,
Robert-Mayer-Strasse 10, 60325 Frankfurt am Main, Germany,
{\tt neuenkir@math.uni-frankfurt.de}
}, 
I. Nourdin\footnote{
Laboratoire de Probabilit\'es et Mod\`eles Al\'eatoires, Universit{\'e} Pierre et Marie Curie,
Bo{\^\i}te courrier 188, 4 Place Jussieu, 75252 Paris Cedex 5, France,
{\tt nourdin@ccr.jussieu.fr}
}, 
A. R\"o{\ss}ler\footnote{
Technical University of Darmstadt, Department of Mathematics,
Schlossgartenstrasse 7, 64289 Darmstadt, Germany,
{\tt roessler@mathematik.tu-darmstadt.de}
}$\,\,$ 
and S. Tindel\footnote{
Institut {\'E}lie Cartan Nancy, B.P. 239,
54506 Vand{\oe}uvre-l{\`e}s-Nancy Cedex, France,
{\tt tindel@iecn.u-nancy.fr}
}
}
\maketitle

\begin{abstract}
In this paper we consider a $n$-dimensional stochastic differential
equation driven by a fractional Brownian motion with Hurst parameter
$H>1/3$. After solving this equation in a rather elementary way,
following the approach of \cite{Gu}, we show how to obtain an expansion
for $E[f(X_t)]$ in terms of $t$, where $X$ denotes the solution to 
the SDE and $f:\R^n\to\R$ is a regular function. 
With respect to \cite{BC}, where the same kind of problem 
is considered, we try an improvement in three different directions:
we are able to take a drift into account in the equation,
we parametrize our expansion with trees (which makes it easier to use),
and we obtain a sharp control of the remainder.
\end{abstract}

\vspace{1cm}


\vspace{0.3cm}


\section{Introduction}
In this article, we  study the stochastic differential equation (SDE in short)
\begin{equation}\label{eqintro}
X_t^a=a+  \int_0^t \sigma(X_s^a)dB_s+ \int_0^t b(X_s^a)ds,\quad t\in [0,T],
\end{equation}
where $B$ is a $d$-dimensional fractional Brownian motion (fBm in short) of Hurst index $H>1/3$,
$a\in\R^n$ is a non-random initial value and $\sigma:\R^n\rightarrow \mathcal{L}^{d,n}$ and $b:\R^n\rightarrow
\R^n$ are smooth functions.

\vspace{0.2cm}

There are essentially two ways to give a sense to equation (\ref{eqintro}):
\begin{enumerate}
\item {\it Pathwise (Stratonovich) setting}.
When $H>1/2$ it is now well-known that we can use the Young integral for integration with respect to fBm and,
with this choice, we have existence and uniqueness of  the solution for equation  (\ref{eqintro}) in the class
of processes having $\alpha$-H\"older continuous paths with
$1-H<\alpha<H$,
see e.g.
\cite{ruz}.
When $1/4<H<1/2$, it is still possible to give a sense to (\ref{eqintro}),
using the rough path theory, which was
initiated by Lyons \cite{LyonsBook,Lyons} and applied to the fBm case by Coutin and Qian \cite{CQ}. In this setting,
we also have existence and uniqueness in an appropriate class of processes.
Remark moreover that, by using a generalization of the symmetric Russo-Vallois integral
(namely the Newton-C\^otes integral corrected by a L\'evy area) we can obtain
existence and uniqueness for (\ref{eqintro}) for  any $H\in (0,1)$, but  only in
dimension $n=d=1$, see \cite{nourdin-simon}.
\item {\it Skorohod setting}.
Skorohod stochastic equations, i.e.,  the integral with respect to
fBm  in (\ref{eqintro}) is understood in the Skorohod sense, are
much more difficult to be solved. Indeed, until now, essentially
only equations in which the noise enters linearly have been
considered, see e.g., \cite{nourdin-tudor}. The difficulty with
equations which are driven  non-linearly by fBm is notorious: the Picard iteration technique involves
Malliavin derivatives in such a way that the equations for
estimating these derivatives cannot be closed.
\end{enumerate}

In the current paper, we will solve (\ref{eqintro}) by means of a variant
of the rough path theory
introduced by Gubinelli in \cite{Gu}.
It is based on an algebraic structure, which
turns out to be useful for computational purposes, but has also its own
interest, and is in fact a nice alternative to the now
classical theory of rough paths initiated by Lyons \cite{LyonsBook,Lyons}.
Although  SDEs of the type (\ref{eqintro}) have already been studied
in \cite{Gu}, we include in this present paper a detailed
review of the algebraic integration tools for several reasons.
First of all, we want to show that this theory can simplify
some aspects of the analysis of fractional equations, and we wish to
give a self-contained study of these objects to illustrate this point.
Moreover, the analysis of  stochastic partial differential equations  in \cite{GT} has lead
to some clarifications with respect to \cite{Gu}, which may be
worth presenting in the simpler finite-dimensional context.
In particular, our computations will heavily rely on an It\^o-type
formula for the so-called weakly controlled processes, which is not included
in \cite{Gu}, and which will be proved here in detail.

As an  application of this theory of integration we study the asymptotic
development with respect to  $t$ of the quantity $P_tf(a)$ 
defined by
\begin{equation}\label{pt}
P_t f(a)={\rm E}(f(X_t^a)),\quad t\in[0,T],\,a\in\R^n,\,f\in {\rm C}^\infty(\R^n;\R),
\end{equation}
where $X^a$ is the solution of (\ref{eqintro}). In the case $H=1/2$,
the Taylor expansion of the semi-group $P_t$ 
is well studied, see, e.g. \cite{Roessler, WP}. Recently, Baudoin and Coutin \cite{BC} studied the asymptotic behaviour in the case $H \neq 1/2$. In this article, we extend their result in several ways:

\begin{enumerate}

\item In \cite{BC}, the authors  considered the particular case $b\equiv 0$.
Consequently, their formula contains only powers of $t$ of the form $t^{nH}$ with
$n\in\N$. Due to the drift part, we obtain  a more complicated expression containing
powers of the type $t^{nH+m}$ with $n,m\in\N$.
\item 
In the current article, we use rooted trees
in order to obtain a nice representation of our formula.
See also \cite{Roessler} for the case $H=1/2$, and \cite{Gu2}
for an application of the tree expansion to the resolution of
stochastic equations.
\item 
In the case where $H>1/2$,  we obtain a series  expansion (\ref{ptfa}) of the operator $P_{t}$, 
which is not only valid for small times as in \cite{BC},
but for any fixed time $t\ge 0$.
\end{enumerate}

Moreover, let us note that in \cite{BC}, the authors used the rough paths theory of Lyons \cite{CQ,LyonsBook,Lyons}
in order to give a sense
to (\ref{eqintro}). Here, as already mentioned, we use the integration theory initiated by Gubinelli \cite{Gu}, 
which allows a self-contained and hopefully a little
simpler version of the essential results contained in the usual theory
of integration of rough signals.

There are several reasons which  motivate the study of the family of
operators $(P_t,\,t\ge 0)$. For instance, the knowledge of $P_t f(a)$ for a sufficiently large
class of functions $f$ characterizes the law of the random variable $X_t^a$.  Moreover, the knowledge of $P_t f(a)$ helps, e.g., also in finding  good sample designs for the reconstruction of fractional diffusions, see \cite{N_in_prep}.

The paper is organized as follows. In Section 2, we state the two
main results of this paper. In section 3,  the basic setup of
\cite{Gu}  with the aim of having a self-contained introduction to
the topic is recalled. In section 4, we recall some facts on 
  the Malliavin calculus for fractional Brownian motion  and some properties of stochastic 
differential equations driven by a fractional  Brownian motion with Hurst parameter $H>1/2$.
Finally, we give the missing proofs in section 5.

\section{Main results}
Before getting into a detailed description of the results
contained in this article,
let us first recall the main properties of a fractional Brownian motion (fBm in short).
A $d$-dimensional fBm with
Hurst parameter $H$ is a
centered Gaussian process, which can be written as
\begin{equation}\label{def:fbm}
B=\lcl B_t=(B_t^1,\ldots,B_t^d);\, t\ge 0 \rcl,
\end{equation}
where $B^1,\ldots,B^d$ are $d$ independent one-dimensional
fBm, i.e., each $B^{i}$ is a centered Gaussian process with continuous sample paths and  covariance function
\begin{equation}\label{rhts}
R_H(t,s)=
\frac12 \lp s^{2H}+t^{2H}-|t-s|^{2H} \rp
\end{equation}
for $i=1, \ldots, d$.
Recall also that  $B^i$ can be represented
in the following way: there exists a standard Brownian motion
$W^i$ such that we have
$$
B_t^i=\iot K(t,s) \,dW^i_s,
$$
 for any $t\ge 0$,
where $K$ is the kernel given by
$$
K(t,s)=c_H
\lc (t-s)^{H-1/2}
+\lp \frac12 -H \rp \ist (u-s)^{H-3/2}
\lp 1- \lp  \frac{s}{u}\rp^{1/2-H} \rp  du \rc \1_{[0,t)}(s),
$$
for a constant $c_H$ which can be expressed in terms of the Gamma function.
Moreover, the fBm verifies the following  two important properties:
\begin{equation}\label{scaling}
\mbox{(scaling)}\quad\mbox{For any }c>0,\,B^{(c)}=c^H B_{\cdot/c}\mbox{ is a fBm},
\end{equation}
\begin{equation}\label{stationarity}
\mbox{(stationarity)}\quad\mbox{For any }h>0,\,B_{\cdot+h}-B_h\mbox{ is a fBm}.
\end{equation}

\bigskip

\subsection{Existence and uniqueness of the solution of  fractional SDEs}

As mentioned in the introduction, we will use for the integration  with respect to fBm  the integration theory developed by Gubinelli \cite{Gu},
on which we try to give here a simplified overview. To this purpose,
will 
denote by $\cl^{d,n}$ the space of linear operators from $\R^d$
to $\R^n$, i.e., the space of matrices of $\R^{n\times d}$. 

The results for fBm we will obtain in section 4 can be summarized as follows:
\begin{theorem}\label{sum-fbm}
Let $B$ be a $d$-dimensional fractional Brownian motion with
Hurst parameter $H>1/3$ and $a\in\R^n$.
Let $b:\R^n\to\R^n$ and $\si:\R^n\to\cl^{d,n}$  be twice continuously differentiable and assume moreover that $\sigma$ and $b$ are bounded
together with their derivatives.
Then the stochastic differential equation
\begin{equation}\label{eds-frac}
X_t^a=a+\iot \si(X_s^a)\, dB_s + \iot b(X_s^a)\, ds, \quad\mbox{
for }\quad t\in\ott,
\end{equation}
admits a unique solution in
$\cq_{\ka,a}(\R^n)$ (see Definition \ref{def39} below) for any $\ka<H$
such that $2\ka+H>1$, where the integral $\iot \si(X^a_s)\, dB_s$
has to be understood in the pathwise sense of Proposition \ref{intg:mdx}.
Moreover, if $f\in C^2(\R^n;\R)$ is bounded together with its
derivatives, then $f(X_t^a)$ can be decomposed as
\begin{equation}\label{ito-frac}
f(X_t^a)= f(a)+ \iot \nabla f(X_s^a)\, b(X_s^a)\, ds +\iot \nabla
f(X_s^a)\, \si(X_s^a)\, dB_s,
\end{equation}
for $t\in\ott$.
\end{theorem}

It is important to note that one of the main differences between our approach and the one
developed in \cite{CQ,LyonsBook} is that the latter heavily relies on
the almost sure approximation of $B$ by a sequence $\{B^n;\, n\ge 1\}$ of
piecewise linear $C^1$-processes, while in our setting this discretization
procedure is only present  for the construction
of the so-called fundamental map $\laa$ (see Proposition \ref{prop:Lambda} below).

\subsection{Rooted trees and their application to the expansion of $P_t$}
To state the next main results we need to recall some properties
of stochastic rooted trees, which have been introduced in \cite{Roessler}
in the case of standard Brownian motion.

\subsubsection{Recalls on rooted trees}
\begin{definition} \label{Def:rooted-S-trees:Wm}
        A {\emph{monotonically labelled S-tree (stochastic
        tree)}} $\textbf{t}$
        with $l=l(\textbf{t}) \in \mathbb{N}$ nodes is a
        pair of maps $\textbf{t}=(\textbf{t}',\textbf{t}'')$
        \begin{equation*}
            \begin{split}
                \textbf{t}' & : \{2, \ldots, l\} \longrightarrow \{1, \ldots,
                l-1\} \\
                \textbf{t}'' & : \{1, \ldots, l\} \longrightarrow
                \mathcal{A}
            \end{split}
        \end{equation*}
        with $\mathcal{A} = \{\gamma, \tau_0, \tau_{j_k}, k \in \mathbb{N}\}$
    where $j_k$ is a variable index with $j_k \in \{1, \ldots, d\}$,
    such that $\textbf{t}'(i) <
        i$, $\textbf{t}''(1)=\gamma$ and $\textbf{t}''(i) \in
        \mathcal{A} \setminus \{\gamma\}$ for $i=2, \ldots, l$. Let $LTS$ denote the
        set of all monotonically labelled S-trees.
\end{definition}
\noindent We use the following notation:
\begin{equation*}
    \begin{split}
        d(\textbf{t}) &= |\{i : \textbf{t}''(i)=\tau_0\}| \\
        s(\textbf{t}) &= |\{i : \textbf{t}''(i)=\tau_{j_k}, j_k \neq 0\}|=l({\bf t})-d({\bf t})-1 \\
        \rho(\textbf{t}) &= H \, s(\textbf{t}) + d(\textbf{t})
    \end{split}
\end{equation*}
with $\rho(\gamma)=0$.

In the following we denote by $LTS(S) \subset LTS$,
where $(S)$ stands for Stratonovich,
the subset
\begin{equation}
    LTS(S) = \{ \textbf{t} \in LTS : s(\textbf{t}) = 2k, \quad
    k \in \mathbb{N}_0 \}
\end{equation}
with $\mathbb{N}_0 = \mathbb{N} \cup \{0\}$ containing all trees
having an even number of stochastic nodes.

\begin{figure}[htb]
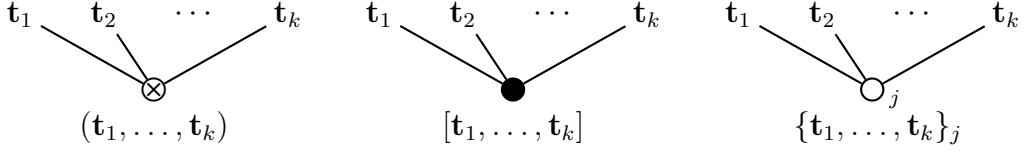

\begin{center}
\begin{tabular}{ccccc}
    \pstree[treemode=U, dotstyle=otimes, dotsize=3.2mm, levelsep=0.1cm, radius=1.6mm, treefit=loose]
    {\Tn}{
    \pstree[treemode=U, dotstyle=otimes, dotsize=3.2mm, levelsep=1cm, radius=1.6mm, treefit=loose, nodesepB=1mm]
    {\Tdot} {\Tr{$\textbf{t}_1$} \Tr{$\textbf{t}_2$} \Tr[edge=none]{$\cdots$} \Tr{$\textbf{t}_k$}}
    }
    & &
    \pstree[treemode=U, dotstyle=otimes, dotsize=3.2mm, levelsep=0.1cm, radius=1.6mm, treefit=loose]
    {\Tn}{
    \pstree[treemode=U, dotstyle=otimes, dotsize=3.2mm, levelsep=1cm, radius=1.6mm, treefit=loose, nodesepB=1mm]
    {\TC*} {\Tr{$\textbf{t}_1$} \Tr{$\textbf{t}_2$} \Tr[edge=none]{$\cdots$} \Tr{$\textbf{t}_k$}}
    }
    & &
    \pstree[treemode=U, dotstyle=otimes, dotsize=3.2mm, levelsep=0.1cm, radius=1.6mm, treefit=loose]
    {\Tn}{
    \pstree[treemode=U, dotstyle=otimes, dotsize=3.2mm, levelsep=1cm, radius=1.6mm, treefit=loose, nodesepB=1mm]
    {\TC~[tnpos=r]{}~[tnpos=r]{$\!\!\!\!_{j}$}} {\Tr{$\textbf{t}_1$} \Tr{$\textbf{t}_2$} \Tr[edge=none]{$\cdots$} \Tr{$\textbf{t}_k$}}
    }
    \\
    $(\textbf{t}_1, \ldots, \textbf{t}_k)$ & &
    $[\textbf{t}_1, \ldots, \textbf{t}_k]$ & & $\,\, \{\textbf{t}_1, \ldots, \textbf{t}_k\}_j$
\end{tabular}
\caption{Writing a coloured S-tree with brackets.}
\label{St-tree-bracket-together}
\end{center}
\end{figure}

Every monotonically labelled S-tree $\textbf{t}$ can be
represented as a graph, whose nodes are elements of $\{1, \ldots,
l(\textbf{t})\}$ and whose arcs are the pairs $(\textbf{t}'(i),i)$
for $i=2, \ldots, l(\textbf{t})$. Here, $\textbf{t}'$ defines a
father son relation between the nodes, i.e., \ $\textbf{t}'(i)$ is
the father of the son $i$.
Further, $\gamma = \,\,$~\pstree[treemode=U, dotstyle=otimes,
dotsize=3.2mm, levelsep=0.1cm, radius=1.6mm, treefit=loose]
    {\Tn}{
    \pstree[treemode=U, dotstyle=otimes, dotsize=3.2mm, levelsep=0cm, radius=1.6mm, treefit=loose]
    {\Tdot~[tnpos=r]{ }} {}
    }
denotes the root, $\tau_0 = $~\pstree[treemode=U, dotstyle=otimes,
dotsize=3.2mm, levelsep=0.1cm, radius=1.6mm, treefit=loose]
    {\Tn}{
    \pstree[treemode=U, dotstyle=otimes, dotsize=3.2mm, levelsep=0cm, radius=1.6mm, treefit=loose]
    {\TC*~[tnpos=r]{}} {}
    }
is a deterministic node and $\tau_{j_k} = $~\pstree[treemode=U,
dotstyle=otimes, dotsize=3.2mm, levelsep=0.1cm, radius=1.6mm,
treefit=loose]
    {\Tn}{
    \pstree[treemode=U, dotstyle=otimes, dotsize=3.2mm, levelsep=0cm, radius=1.6mm, treefit=loose]
    {\TC~[tnpos=r]{$\!\!_{j_k}$}} {}
    }
a stochastic node. Here, we have to point out that each tree
$\textbf{t} \in LTS$ depends on the variable indices $j_1, \ldots,
j_{s(\textbf{t})} \in \{1, \ldots, d\}^{s(\textbf{t})}$, although
this is not mentioned explicitly if we shortly write $\textbf{t}$
for the tree.
\begin{figure}[htbp]
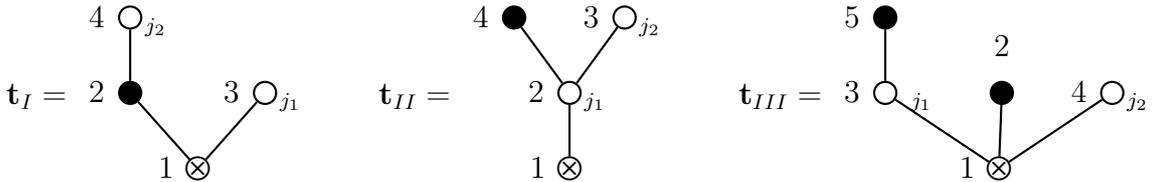

\begin{center}
\begin{tabular}{ccccc}
    $\textbf{t}_I = \begin{array}{c}
    \text{\pstree[treemode=U, dotstyle=otimes, dotsize=3.2mm, levelsep=0.1cm, radius=1.6mm, treefit=loose]
    {\Tn}{
    \pstree[treemode=U, dotstyle=otimes, dotsize=3.2mm, levelsep=1cm, radius=1.6mm, treefit=loose]
    {\Tdot~[tnpos=l]{1 }}
    {\pstree{\TC*~[tnpos=l]{2}}{\TC~[tnpos=l]{4}~[tnpos=r]{$\!\!_{j_2}$}}
    \TC~[tnpos=l]{3}~[tnpos=r]{$\!\!_{j_1}$}}
    }}
    \end{array}$
    & \quad  &
    $\textbf{t}_{II} = \begin{array}{c}
    \text{\pstree[treemode=U, dotstyle=otimes, dotsize=3.2mm, levelsep=0.1cm, radius=1.6mm, treefit=loose]
    {\Tn}{
    \pstree[treemode=U, dotstyle=otimes, dotsize=3.2mm, levelsep=1cm, radius=1.6mm, treefit=loose]
    {\Tdot~[tnpos=l]{1 }} {\pstree{\TC~[tnpos=l]{2}~[tnpos=r]{$\!\!_{j_1}$}}{\TC*~[tnpos=l]{4}
    \TC~[tnpos=l]{3}~[tnpos=r]{$\!\!_{j_2}$}}}
    }}
    \end{array}$
    & \quad  &
    $\textbf{t}_{III} = \begin{array}{c}
    \text{\pstree[treemode=U, dotstyle=otimes, dotsize=3.2mm, levelsep=0.1cm, radius=1.6mm, treefit=loose]
    {\Tn}{
    \pstree[treemode=U, dotstyle=otimes, dotsize=3.2mm, levelsep=1cm, radius=1.6mm, treefit=loose]
    {\Tdot~[tnpos=l]{1 }} {\pstree{\TC~[tnpos=l]{3}~[tnpos=r]{$_{j_1}$}}{\TC*~[tnpos=l]{5}} \TC*~[tnpos=a]{2} \TC~[tnpos=l]{4}~[tnpos=r]{$\!\!_{j_2}$}}
    }}
    \end{array}$
\end{tabular}
\caption{Some monotonically labelled trees in $LTS$.}
\label{St-S-tree-examples-tI+tII:Wm}
\end{center}
\end{figure}
\begin{definition}
    If $\textbf{t}_1, \ldots, \textbf{t}_k$ are coloured trees then we denote by
    \begin{equation*}
        (\textbf{t}_1, \ldots, \textbf{t}_k), \,\,\,\,\,
        [\textbf{t}_1, \ldots, \textbf{t}_k] \,\,\,\,\, \text{ and }
        \,\,\,\,\, \{\textbf{t}_1, \ldots, \textbf{t}_k \}_{j}
    \end{equation*}
    the tree in which $\textbf{t}_1, \ldots, \textbf{t}_k$ are each joined by a
    single branch to $\,\,$
    \pstree[treemode=U, dotstyle=otimes, dotsize=3.2mm, levelsep=0.1cm, radius=1.6mm, treefit=loose]
    {\Tn}{
    \pstree[treemode=U, dotstyle=otimes, dotsize=3.2mm, levelsep=0cm, radius=1.6mm, treefit=loose]
    {\Tdot} {}
    }
    $\,$,
    \pstree[treemode=U, dotstyle=otimes, dotsize=3.2mm, levelsep=0.1cm, radius=1.6mm, treefit=loose]
    {\Tn}{
    \pstree[treemode=U, dotstyle=otimes, dotsize=3.2mm, levelsep=0cm, radius=1.6mm, treefit=loose]
    {\TC*} {}
    }
    $\,\,$and
    \pstree[treemode=U, dotstyle=otimes, dotsize=3.2mm, levelsep=0.1cm, radius=1.6mm, treefit=loose]
    {\Tn}{
    \pstree[treemode=U, dotstyle=otimes, dotsize=3.2mm, levelsep=0cm, radius=1.6mm, treefit=loose]
    {\TC~[tnpos=r]{}~[tnpos=r]{$\!\!\!\!_{j}$}} {}
    },
    respectively (see also Figure~\ref{St-tree-bracket-together}).
\end{definition}
Therefore proceeding recursively, for the two examples
$\textbf{t}_I$ and $\textbf{t}_{II}$ in
Figure~\ref{St-S-tree-examples-tI+tII:Wm} we obtain
    $\textbf{t}_I = ([
    \text{\pstree[treemode=U, dotstyle=otimes, dotsize=3.2mm, levelsep=0.1cm, radius=1.6mm, treefit=loose]
    {\Tn}{
    \pstree[treemode=U, dotstyle=otimes, dotsize=3.2mm, levelsep=0cm, radius=1.6mm, treefit=loose]
    {\TC~[tnpos=r]{$\!^4_{j_2}$}} {}
    }} ]^2,
    \text{\pstree[treemode=U, dotstyle=otimes, dotsize=3.2mm, levelsep=0.1cm, radius=1.6mm, treefit=loose]
    {\Tn}{
    \pstree[treemode=U, dotstyle=otimes, dotsize=3.2mm, levelsep=0cm, radius=1.6mm, treefit=loose]
    {\TC~[tnpos=r]{$\!^3_{j_1}$}} {}
    }} )^1 = ([\tau_{j_2}^4]^2 , \tau_{j_1}^3)^1$
    and
    $\textbf{t}_{II} = (\{
    \text{\pstree[treemode=U, dotstyle=otimes, dotsize=3.2mm, levelsep=0.1cm, radius=1.6mm, treefit=loose]
    {\Tn}{
    \pstree[treemode=U, dotstyle=otimes, dotsize=3.2mm, levelsep=0cm, radius=1.6mm, treefit=loose]
    {\TC*~[tnpos=r]{$\!^4$}} {}
    }},
    \text{\pstree[treemode=U, dotstyle=otimes, dotsize=3.2mm, levelsep=0.1cm, radius=1.6mm, treefit=loose]
    {\Tn}{
    \pstree[treemode=U, dotstyle=otimes, dotsize=3.2mm, levelsep=0cm, radius=1.6mm, treefit=loose]
    {\TC~[tnpos=r]{$\!^3_{j_2}$}} {}
    }}
    \}_{j_1}^2)^1 = (\{ \tau_0^4, \tau_{j_2}^3 \}_{j_1}^2)^1$. 

\vspace{0.2cm}

For every rooted tree $\textbf{t} \in LTS$, there exists a
corresponding elementary differential which is a direct
generalization of the differential in the deterministic case, see
also \cite{Roessler}. The elementary differential is defined
recursively for some $x \in \mathbb{R}^n$ by
\begin{equation*}
    F(\gamma)(x) = f(x), \qquad
    F(\tau_0)(x) = b(x), \qquad
    F(\tau_j)(x) = \sigma^j(x),
\end{equation*}
for single nodes and by
\begin{equation} \label{St-elementary-differential-F:Wm}
    F(\textbf{t})(x) =
    \begin{cases}
    f^{(k)}(x) \cdot (F(\textbf{t}_1)(x), \ldots, F(\textbf{t}_k)(x)) &
    \text{for } \textbf{t}=(\textbf{t}_1, \ldots, \textbf{t}_k) \\
    b^{(k)}(x) \cdot (F(\textbf{t}_1)(x), \ldots,
    F(\textbf{t}_k)(x)) & \text{for } \textbf{t}=[\textbf{t}_1, \ldots, \textbf{t}_k] \\
    {\sigma^j}^{(k)}(x) \cdot (F(\textbf{t}_1)(x),
    \ldots, F(\textbf{t}_k)(x)) & \text{for } \textbf{t}=\{\textbf{t}_1, \ldots,
    \textbf{t}_k\}_j
    \end{cases}
\end{equation}
for a tree $\textbf{t}$ with more than one node and with $\sigma^j
= (\sigma^{i,j})_{1 \leq i \leq n}$ denoting the $j$th column of
the diffusion matrix $\sigma$. Here $f^{(k)}$, $b^{(k)}$ and
${\sigma^j}^{(k)}$ define a symmetric $k$-linear differential
operator, and one can choose the sequence of labelled S-trees
$\textbf{t}_1, \ldots, \textbf{t}_k$ in an arbitrary order. For
example, the $I$th component of $b^{(k)} \cdot (F(\textbf{t}_1),
\ldots, F(\textbf{t}_k))$ can be written as
\begin{equation*}
    ( b^{(k)} \cdot (F(\textbf{t}_1), \ldots, F(\textbf{t}_k)) )^I
    = \sum_{J_1, \ldots, J_k=1}^n \frac{\partial^k
    b^I}{\partial x^{J_1} \ldots \partial x^{J_k}} \,
    (F^{J_1}(\textbf{t}_1), \ldots, F^{J_k}(\textbf{t}_k)),
\end{equation*}
where the components of vectors are denoted by superscript
indices, which are chosen as capitals.
%
%
As a result of this we get for $\textbf{t}_I$ and
$\textbf{t}_{II}$ the elementary differentials
\begin{equation*}
    \begin{split}
    &F(\textbf{t}_I) = f'' (b'(\sigma^{j_2}), \sigma^{j_1}) = \sum_{J_1,J_2=1}^n
    \frac{\partial^2 f}{\partial x^{J_1} \partial x^{J_2}}
    \left( \sum_{K_1=1}^n \frac{\partial b^{J_1}}{\partial x^{K_1}}
    \, \sigma^{K_1,j_2} \cdot \sigma^{J_2,j_1} \right), \\
    &F(\textbf{t}_{II}) = f' ({\sigma^{j_1}}'' (b, \sigma^{j_2})) = \sum_{J_1=1}^n
    \frac{\partial f}{\partial x^{J_1}} \left( \sum_{K_1, K_2 =1}^n
    \frac{\partial^2 \sigma^{J_1,j_1}}{\partial x^{K_1} \partial
    x^{K_2}} \, b^{K_1} \cdot \sigma^{K_2,j_2} \right).
    \end{split}
\end{equation*}
Next, we assign recursively to every $\textbf{t} \in LTS$ a
multiple stochastic integral
by
\begin{equation} \label{Def-tree-stoch-int}
    \mathcal{I}_{\textbf{t}}(g(X_s^a))_{t_0,t} = \begin{cases}
        \displaystyle
        g(X_t^a)
        & \text{ if } \textbf{t}''(l(\textbf{t})) = \gamma \\
        \displaystyle
        \int_{t_0}^t \mathcal{I}_{\textbf{t}-}(g(X_u^a))_{t_0,s} \, dB_s^j
        & \text{ if } \textbf{t}''(l(\textbf{t})) = \tau_j
        \end{cases}
\end{equation}
for $0 \leq j \leq d$ with $dB_s^0=ds$. Here, $\textbf{t}-$ denotes the tree which is obtained from
$\textbf{t}$ by removing the last node with label $l(\textbf{t})$. 
\subsubsection{Expansion of $P_t$ with respect to time $t$}

We will denote by $C^{\infty}_{b}(\R^{n};\R)$ the space  of all infinitely differentiable 
functions $g \in C^{\infty}(\R^{n};\R)$, which are bounded together with their derivatives. 
Moreover set $\mathcal{A}_m=\{0,1,\ldots,d\}^m$ for $m\in\N$, and define the differential operators
$\mathcal{D}^0$ and
$\mathcal{D}^j$ as
\begin{equation}
    \mathcal{D}^0 = \sum_{k=1}^n b^k \, \frac{\partial}{\partial x^k} \qquad \qquad
    \text{ and } \qquad \qquad \mathcal{D}^j = \sum_{k=1}^n \sigma^{k,j} \,
    \frac{\partial}{\partial x^k}
\end{equation}
for $j=1, \ldots, d$. Finally, set 
$\mathcal{D}^{\alpha} = \mathcal{D}^{\alpha_1} \ldots
\mathcal{D}^{\alpha_m}$ for a multi-index $\alpha \in \mathcal{A}_{m}$. 

Recall that the family of operators $(P_t,\,t\in [0,T])$ has been
defined by (\ref{pt}).
To give the expansion of $P_{t}$ we will
use the following assumptions on the function $f: \R^{n}
\rightarrow \R$, the drift vector $b=(b^{i})_{i=1, \ldots, n} $
and the diffusion matrix $\sigma= (\sigma^{i,j})_{i=1, \ldots, n,
\,  j=1, \ldots d}$:

\medskip

\begin{itemize}
\item[(A)] We have $f, b^{i}, \sigma^{i,j} \in C^{\infty}_{b}(\R^{n}; \R)$ for $i=1, \ldots, n$, $j=1, \ldots, d$.
\end{itemize}

The following theorem gives the expression of the expansion of
$P_t$ with respect to $t$:

\begin{theorem}\label{thmptfa}
\begin{enumerate}
\item If $H>1/3$ and assumption (A) is satisfied, we have that
\begin{equation}\label{ptfasmall}
     {\rm P}_t f(a) = \sum_{\substack{\textbf{t} \in LTS(S) \\
        l(\textbf{t}) \leq m+1}}
        \sum_{j_1, \ldots, j_{s(\textbf{t})}=1}^d
        F(\textbf{t})(a) \,\, {\rm E}(\mathcal{I}_{\textbf{t}}(1)_{0,1})
        \,\, t^{\rho(\textbf{t})} + O(t^{(m+1)H}),\mbox{ as }t\rightarrow 0
\end{equation}
for any $m \in \N_0$.
\item Let $H>1/2$ and  assumption (A) be satisfied. 
Moreover, assume that there exist $M \in \N$ and  constants $K>0$, $\kappa \in [0,1/2)$ such 
that
\begin{equation}\label{add_ass} 
\sup_{\alpha\in\mathcal{A}_m}\sup_{x\in\R^n}|{\mathcal D}^\alpha f(x)|
+\max_{i=1,\ldots,n}\sup_{\alpha\in\mathcal{A}_m}\sup_{x\in\R^n}\left|\frac{\partial}
{\partial x_i}{\mathcal D}^\alpha f(x)\right|
\leq K^{m} (m!)^{\kappa}
\end{equation}  
  for all  $m \geq M$. 
Then we have 
\begin{equation}\label{ptfa}
        {\rm P}_t f(a) = \sum_{\textbf{t} \in LTS(S)}
        \sum_{j_1, \ldots, j_{s(\textbf{t})}=1}^d
        F(\textbf{t})(a) \,\, {\rm E}(\mathcal{I}_{\textbf{t}}(1)_{0,1})
        \,\, t^{\rho(\textbf{t})},\quad t\in [0,T].
\end{equation}
\end{enumerate}
\end{theorem}

\bigskip

\begin{remark} 
{\rm 
{\bf (1)} Here, note that each tree $\textbf{t} \in LTS(S)$ comprehends the
variable indices $j_1, \ldots, j_{s(\textbf{t})}$ which can take
the values $1, \ldots, d$ although these variables are not
mentioned explicitly by writing shortly $\textbf{t}$ for the whole
tree. The variables $j_1, \ldots, j_{s(\textbf{t})}$ correspond to
the components of the driving fractional Brownian motion and
appear in the second sum in the formulas (\ref{ptfasmall}) and
(\ref{ptfa}) as well as in each tree $\textbf{t}$ of the summands.

\vspace{0.3cm}

\noindent
{\bf (2)}
In the case where $H>1/2$, the boundedness of the coefficients is not needed 
for existence and uniqueness of the solution, see \cite{NR}.

\vspace{0.3cm}

\noindent
{\bf (3)}
Although the additional assumption (\ref{add_ass}) seems to   
be quite restrictive, it is however natural in a certain sense.  
Indeed, consider the trivial one-dimensional equation
$$ dX_{t}^{a} = dB_{t}, \quad t \in [0,T], \qquad X_0=a$$
for $H>1/2$.
Then we have clearly
$$ X_{t}^{a}=a+B_{t}, \qquad t \in [0,T].$$
By the first part of Theorem \ref{thmptfa}  and Proposition \ref{mom_12} we have for  this equation the expansion 
$$ {\rm E } f(X_{t}^{a}) = f(a) + \sum_{k=1}^{m} c_{k}f^{(k)}(a)  t^{Hk} + O(t^{H(m+1)})$$
with 
\begin{eqnarray*} c_{k}=  \frac{{\rm E}(B_{1})^{k}}{k!} = \left \{ \begin{array}{cc} 0, & \textrm{if } k \textrm{ is odd},  \\  \frac{1}{2^{k/2} (k/2)!} & \textrm{ else.} \end{array} \right. \end{eqnarray*}
Hence the series 
$$  \sum_{k=1}^{\infty} c_{k}f^{(k)}(a)  t^{Hk}$$
converges absolutely for any $t\in [0,T]$,
for instance if there exists constants $K>0$ and $\kappa \in [0,1/2)$ such that
$$ |f^{(m)}(a)| \leq K^{m} (m!)^{\kappa}$$
for all $m \in \N$.
But if, e.g., we have that
$$ \liminf_{m \rightarrow \infty} \frac{|f^{(m)}(a)|}{ (m!)^{1/2+\varepsilon}} >0 $$ with $\varepsilon >0$, then we have
$$ \sum_{k=1}^{\infty} c_{k}|f^{(k)}(a)|  t^{Hk} = + \infty  $$ for any $t \in (0,T]$.  
Thus the  condition (\ref{add_ass}) we require for the control of the remainder is   quite natural, since the coefficients of the expansion have to satisfy a similar condition, as illustrated in this example. 
Furthermore similar growth conditions on the remainder or the coefficients are also 
usual in the case $H=1/2$, i.e., for the asymptotic expansion  of It\^{o} stochastic differential 
equations  respectively their functionals. Compare, e.g., \cite{BA} and chapter 5 in \cite{KP}. 

\vspace{0.3cm}

\noindent
{\bf (4)}
In order to solve equation (\ref{eds-frac}) and to bound the 
Malliavin derivative of the solution in the case $H>1/2$, we need  only
a boundedness condition on the first two derivatives of $b$ and $\si$.
To avoid too many technicalities, 
we have assumed in (\ref{add_ass})
that all derivatives $\mathcal D^\al f$ are  uniformly bounded in $x\in\R^n$.
However, thanks to Proposition \ref{prop_X} part (b), this condition
could be relaxed, and we could allow a bound of the form
$$
\sup_{\alpha\in\mathcal{A}_m}|{\mathcal D}^\alpha f(x)|
+\max_{i=1,\ldots,n}\sup_{\alpha\in\mathcal{A}_m}\left|\frac{\partial}
{\partial x_i}{\mathcal D}^\alpha f(x)\right|
\leq K^{m} (m!)^{\kappa}\lp 1+\|x\|^q \rp,
$$
for a given $q\ge 0$.
}
\end{remark}

For a better understanding of the previous results and as an
example we consider SDE~(\ref{eqintro}) in the case of $ n \geq
1$, $d \geq 1$ and give an expansion of ${\rm P}_t f(a)$ for
$m=2$. Here, we have to consider the trees with $l(\textbf{t})
\leq 3$ which are $\textbf{t}_{1}=\gamma^1$, $\textbf{t}_{2}
=(\sigma_{j_1}^2)^1$, $\textbf{t}_{3}=(\tau^2)^1$,
$\textbf{t}_{4}=(\sigma_{j_1}^2, \sigma_{j_2}^3)^1$,
$\textbf{t}_{5}=(\{\sigma_{j_2}^3\}_{j_1}^2)^1$,
$\textbf{t}_{6}=([\sigma_{j_1}^3]^2)^1$,
$\textbf{t}_{7}=(\{\tau^3\}_{j_1}^2)^1$,
$\textbf{t}_{8}=(\tau^2,\sigma_{j_1}^3)^1$,
$\textbf{t}_{9}=(\tau^3,\sigma_{j_1}^2)^1$,
$\textbf{t}_{10}=(\tau^2, \tau^3)^1$ and
$\textbf{t}_{11}=([\tau^3]^2)^1$. However, only trees in $LTS(S)$
with an even number of stochastic nodes have to be included since
we have $\E(\mathcal{I}_{\textbf{t}}(1)_{0,1})=0$ for $\textbf{t}
\in LTS \setminus LTS(S)$. Then, we obtain
\begin{align*}
    {\rm P}_t f(a) &= F(\textbf{t}_{1})(a) + F(\textbf{t}_{3})(a) \, {\rm
    E}(\mathcal{I}_{\textbf{t}_{3}}(1)_{0,1}) \, t + \sum_{j_1,j_2=1}^d F(\textbf{t}_{4})(a) \, {\rm
    E}(\mathcal{I}_{\textbf{t}_{4}}(1)_{0,1}) \, t^{2H} \\
    &+ \sum_{j_1,j_2=1}^d F(\textbf{t}_{5})(a) \, {\rm
    E}(\mathcal{I}_{\textbf{t}_{5}}(1)_{0,1}) \, t^{2H}
    + F(\textbf{t}_{10})(a) \, {\rm
    E}(\mathcal{I}_{\textbf{t}_{10}}(1)_{0,1}) \, t^{2} \\
    &+ F(\textbf{t}_{11})(a) \, {\rm
    E}(\mathcal{I}_{\textbf{t}_{11}}(1)_{0,1}) \, t^{2} + O(t^{3H}) \, .
\end{align*}
Applying now (\ref{St-elementary-differential-F:Wm}) and
(\ref{Def-tree-stoch-int}) yields
\begin{align*}
    {\rm P}_t f(a) &= f(a) + f'(b)(a) \, {\rm
    E}(\int_0^1 \, ds) \, t + \sum_{j_1,j_2=1}^d f''(\sigma^{j_1}, \sigma^{j_2})(a) \, {\rm
    E}(\int_0^1 \int_0^s \, dB^{j_1}_{s_1} \, dB^{j_2}_s) \, t^{2H} \\ &+ \sum_{j_1,j_2=1}^d
    f'({\sigma^{j_1}}'(\sigma^{j_2}))(a) \, {\rm
    E}(\int_0^1 \int_0^s \, dB^{j_1}_{s_1} \, dB^{j_2}_s) \, t^{2H} + f''(b,b)(a) \, {\rm
    E}(\int_0^1 \int_0^s \, ds_1 \, ds) \, t^{2} \\ &+ f'(b'(b))(a) \, {\rm
    E}(\int_0^1 \int_0^s \, ds_1 \, ds) \, t^{2} + O(t^{3H}) \, ,
\end{align*}
which finally results in
\begin{align*}
    {\rm P}_t f(a) &= f(a) + \sum_{J_1=1}^n \frac{\partial f}{\partial x^{J_1}}(a) \, b^{J_1}(a) \, t
    +  \frac{1}{2}  \sum_{j=1}^d \sum_{J_1,J_2=1}^n \frac{\partial^2 f}{\partial x^{J_1} \partial x^{J_2}}(a) \,
    \sigma^{J_1,j}(a) \, \sigma^{J_2,j}(a) \, t^{2H} \\
    &+  \frac{1}{2} \sum_{j=1}^d \sum_{J_1,J_2=1}^n \frac{\partial f}{\partial x^{J_1}}(a) \,
    \frac{\partial \sigma^{J_1,j}}{\partial x^{J_2}}(a)
    \, \sigma^{J_2,j}(a) \, t^{2H} \\
    &+  \frac{1}{2} \sum_{J_1,J_2=1}^n \frac{\partial^2 f}{\partial x^{J_1} \partial
    x^{J_2}}(a) \, b^{J_1}(a) \, b^{J_2}(a) \,  t^{2} \\
    &+  \frac{1}{2} \sum_{J_1,J_2=1}^n \frac{\partial f}{\partial x^{J_1}}(a) \, \frac{\partial b^{J_1}}{\partial x^{J_2}}(a)
    \, b^{J_2}(a) \, t^{2} + O(t^{3H}) \, .
\end{align*}


\section{Some elements of algebraic integration}
\label{sec:one-dim}
As  already mentioned  in the introduction, we include
in this present section a detailed
review of the algebraic integration tools contained
mostly in \cite{Gu,GT}.
Moreover, we will give  an  It\^o's type
formula for  so-called weakly controlled processes.

\subsection{Increments}\label{incr}

The extended pathwise integration we will deal with is based on
the notion of 'increments', together with an
elementary operator $\der$ acting on them.
The algebraic structure they generate is described in \cite{Gu,GT},
but here we present  directly the definitions of
interest for us, for sake of conciseness.
First of all,  for an arbitrary real number
$T>0$, a vector space $V$ and an integer $k\ge 1$ we denote by
$\cac_k(V)$ the set of functions $g : [0,T]^{k} \to V$ such
that $g_{t_1 \cdots t_{k}} = 0$
whenever $t_i = t_{i+1}$ for some $i\le k-1$.
Such a function will be called a
\emph{$(k-1)$-increment}, and we will
set $\cac_*(V)=\cup_{k\ge 1}\cac_k(V)$. We can now define the announced
elementary operator $\der$ on $\cac_k(V)$:
\begin{equation}
  \label{eq:coboundary}
\delta : \cac_k(V) \to \cac_{k+1}(V), \qquad
(\delta g)_{t_1 \cdots t_{k+1}} = \sum_{i=1}^{k+1} (-1)^{k-i}
g_{t_1  \cdots \hat t_i \cdots t_{k+1}} ,
\end{equation}
where $\hat t_i$ means that this particular argument is omitted.
A fundamental property of $\der$, which is easily verified,
is that
$\delta \delta = 0$, where $\delta \delta$ is considered as an operator
from $\cac_k(V)$ to $\cac_{k+2}(V)$.
 We will denote $\cz\cac_k(V) = \cac_k(V) \cap \text{Ker}\delta$
and $\cb \cac_k(V) =
\cac_k(V) \cap \text{Im}\delta$.

\vspace{0.3cm}

Some simple examples of actions of $\der$,
which will be the ones we will really use throughout the paper,
 are obtained by letting
$g\in\cac_1$ and $h\in\cac_2$. Then, for any $s,u,t\in\ott$, we have
\begin{equation}
\label{eq:simple_application}
  (\der g)_{st} = g_t - g_s,
\quad\mbox{ and }\quad
(\der h)_{sut} = h_{st}-h_{su}-h_{ut}.
\end{equation}
Furthermore, it is easily checked that
$\cz \cac_{k+1}(V) = \cb \cac_{k}(V)$ for any $k\ge 1$.
In particular, the following basic property holds:
\begin{lemma}\label{exd}
Let $k\ge 1$ and $h\in \cz\cac_{k+1}(V)$. Then there exists a (non unique)
$f\in\cac_{k}(V)$ such that $h=\der f$.
\end{lemma}

\noindent{\it Proof}.
This elementary proof is included in \cite{Gu}, and will be omitted
here. Let us just mention that
$f_{t_1\ldots t_{k}}=h_{0t_1\ldots t_{k}}$ is a possible choice.
\fin

Observe that Lemma \ref{exd} implies that all the elements
$h \in\cac_2(V)$ such that $\der h= 0$ can be written as $h = \der f$
for some (non unique) $f \in \cac_1(V)$. Thus we get a heuristic
interpretation of $\der |_{\cac_2(V)}$:  it measures how much a
given 1-increment  is far from being an  exact increment of a
function, i.e., a finite difference.

\vspace{0.3cm}

Notice that our future discussions will mainly rely on
$k$-increments with $k \le 2$, for which we will use some
analytical assumptions. Namely,
we measure the size of these increments by H\"older norms
defined in the following way: for $f \in \cac_2(V)$ let
$$
\norm{f}_{\mu} =
\sup_{s,t\in\ott}\frac{|f_{st}|}{|t-s|^\mu},
\quad\mbox{and}\quad
\cac_2^\mu(V)=\lcl f \in \cac_2(V);\, \norm{f}_{\mu}<\infty  \rcl.
$$
Obviously, the usual H\"older spaces $\cac_1^\mu(V)$ will be determined
        in the following way: for a continuous function $g\in\cac_1(V)$, we simply set
\begin{equation}\label{def:hnorm-c1}
\|g\|_{\mu}=|\der g|_{\mu},
\end{equation}
and we will say that $g\in\cac_1^\mu(V)$ iff $\|g\|_{\mu}$ is finite.
Notice that $\|\cdot\|_{\mu}$ is only a semi-norm on $\cac_1(V)$,
but we will generally work on spaces of the type
\begin{equation}\label{def:hold-init}
\cac_{1,a}^\mu(V)=
\lcl g:\ott\to V;\, g_0=a,\, \|g\|_{\mu}<\infty \rcl,
\end{equation}
for a given $a\in V$, on which $\|g\|_{\mu}$ thus becomes a norm.
 For $h \in \cac_3(V)$ set in the same way
\begin{eqnarray}
  \label{eq:normOCC2}
  \norm{h}_{\gamma,\rho} &=& \sup_{s,u,t\in\ott}
\frac{|h_{sut}|}{|u-s|^\gamma |t-u|^\rho}\\
\|h\|_\mu &= &
\inf\left \{\sum_i \|h_i\|_{\rho_i,\mu-\rho_i} ;\, h =
 \sum_i h_i,\, 0 < \rho_i < \mu \right\} ,\nonumber
\end{eqnarray}
where the last infimum is taken over all sequences $\{h_i \in \cac_3(V) \}$
such that $h
= \sum_i h_i$ and for all choices of the numbers $\rho_i \in (0,z)$. 
Then  $\|\cdot\|_\mu$ is easily seen to be a norm on $\cac_3(V)$, and we set
$$
\cac_3^\mu(V):=\lcl h\in\cac_3(V);\, \|h\|_\mu<\infty \rcl.
$$
Eventually,
let $\cac_3^{1+}(V) = \cup_{\mu > 1} \cac_3^\mu(V)$,
and remark that the same kind of norms can be considered on the
spaces $\cz \cac_3(V)$, leading to the definition of some spaces
$\cz \cac_3^\mu(V)$ and $\cz \cac_3^{1+}(V)$.

\vspace{0.3cm}

With these notations in mind
the following proposition is a basic result, which  belongs to  the core of
 our approach to pathwise integration. Its proof may be found
in a simple form in \cite{GT}.
\begin{proposition}[The $\Lambda$-map]
\label{prop:Lambda}
There exists a unique linear map $\Lambda: \cz \cac^{1+}_3(V)
\to \cac_2^{1+}(V)$ such that
$$
\delta \Lambda  = \id_{\cz \cac_3^{1+}(V)}
\quad \mbox{ and } \quad \quad
\Lambda  \delta= \id_{\cac_2^{1+}(V)}.
$$
In other words, for any $h\in\cac^{1+}_3(V)$ such that $\der h=0$
there exists a unique $g=\laa(h)\in\cac_2^{1+}(V)$ such that $\der g=h$.
Furthermore, for any $\mu > 1$,
the map $\laa$ is continuous from $\cz \cac^{\mu}_3(V)$
to $\cac_2^{\mu}(V)$ and we have
\begin{equation}\label{ineqla}
\|\Lambda h\|_{\mu} \le \frac{1}{2^\mu-2} \|h\|_{\mu} ,\qquad h \in
\cz \cac^{\mu}_3(V).
\end{equation}
\end{proposition}

\vspace{0.3cm}

We can now give an algorithm for a canonical decomposition of a function
$g\in\cac_2 (V)$, whose increment $\der g$ is smooth enough:
\begin{corollary}
Let $g\in\cac_2 (V)$ such that $\der g\in\cac_3^\mu(V)$ for
$\mu>1$. Then, for an arbitrary $a\in V$,  $g$ can be
decomposed in a unique way as
\begin{equation}\label{dcp:preimage}
g=\der f+ \Lambda \delta g,
\end{equation}
where $f\in\cac_{1,a}(V)$.
\end{corollary}

\noindent{\it Proof}.
This proof is elementary. We include it here in order to
see some simple manipulations of the objects we have introduced
so far.

\vspace{0.3cm}

The existence of the decomposition is due to the following fact:
if $\der g\in\cac_3^\mu(V)$, then it belongs  to $\mbox{Dom}(\laa)$. Thus, let us set $h=\laa\der g$. Then $\der(g-h)=0$, which  means
that $g-h\in\cz\cac_2$, and since $\cz\cac_2=\cb\cac_1$, there exists
an element $f\in\cac_1$ such that $g-h=\der f$. Hence we have
obtained a decomposition of the form (\ref{dcp:preimage}).

\vspace{0.3cm}

As far as the uniqueness of the decomposition is concerned, if $f^1,f^2$ satisfy
(\ref{dcp:preimage}), then $\der f^1=\der f^2$ and hence
they differ only  by a constant. Since $f^1,f^2$ are both supposed
to be elements of $\cac_{1,a}(V)$, where $a$ is a fixed initial
condition, we obtain $f^1=f^2$, which proves our claim.

\fin

\vspace{0.3cm}

Let us mention at this point a first link between the structures
we have introduced so far and the problem of integration of irregular
functions.

\begin{corollary}
\label{cor:integration}
For any 1-increment $g\in\cac_2 (V)$ such that $\der g\in\cac_3^{1+}$,
set
$
\delta f = (\id-\Lambda \delta) g
$.
Then
$$
(\delta f)_{st} = \lim_{|\Pi_{ts}| \to 0} \sum_{i=0}^n g_{t_i\,t_{i+1}},
$$
where the limit is over any partition $\Pi_{ts} = \{t_0=t,\dots,
t_n=s\}$ of $[t,s]$, whose mesh tends to zero. Thus, the
1-increment $\delta f$ is the indefinite integral of the 1-increment $g$.
\end{corollary}

\noindent{\it Proof}.
Just consider the equation $g = \delta f + \Lambda \delta g$ and write
\begin{equation*}
  \begin{split}
S_\Pi & = \sum_{i=0}^n g_{t_i\, t_{i+1}}
=
\sum_{i=0}^n (\delta f)_{t_i\, t_{i+1}}
+
\sum_{i=0}^n (\Lambda \delta g)_{t_i\, t_{i+1}}
\\ & =
(\delta f)_{st}
+
\sum_{i=0}^n (\Lambda \delta g)_{t_i\, t_{i+1}}.
      \end{split}
\end{equation*}
Then observe that, due to the fact that $\Lambda \delta g \in
\cac_3^{1+}(V)$, the last sum converges to zero.
\fin

\subsection{Computations in $\cac_*$}\label{cpss}

Let us specialize now to the case $V=\R^d$ for $d\ge 1$. We will also
denote by $\cl^{d,l}$ the space of linear operators from $\R^d$
to $\R^l$, i.e., the space of matrices of $\R^{l\times d}$
and set $\cac_k\cl^{d,l}=\cac_k(\cl^{d,l})$.
Then $(\cac_*,\delta)$ can be endowed with the following product:
for  $g\in\cac_n\cl^{d,l}$ and $h\in\cac_m(\R^d) $ let  $gh$
be the element of $\cac_{n+m-1}(\R^l)$ defined by
\begin{equation}\label{cvpdt}
(gh)_{t_1,\dots,t_{m+n+1}}=
g_{t_1,\dots,t_{n}} h_{t_{n},\dots,t_{m+n-1}},
\quad
t_1,\dots,t_{m+n-1}\in\ott.
\end{equation}
In this context, we have the following useful properties.

\begin{proposition}\label{difrul}
The following differentiation rules hold true:
\begin{enumerate}
\item
Let $g\in\cac_1\cl^{d,l}$ and $h\in\cac_1(\R^d)$. Then
$gh\in\cac_1(\R^l)$ and
\begin{equation}\label{difrulu}
\der (gh) = \der g\,  h + g\, \der h.
\end{equation}
\item
Let $g\in\cac_1\cl^{d,l}$ and $h\in\cac_2(\R^d)$. Then
$gh\in\cac_2(\R^l)$ and
\begin{equation}
\der (gh) = \der g\, h - g \,\der h.
\end{equation}
\item
Let $g\in\cac_2\cl^{d,l}$ and $h\in\cac_1(\R^d)$. Then
$gh\in\cac_2(\R^l)$ and
\begin{equation}
\der (gh) = \der g\, h  + g \,\der h.
\end{equation}
\end{enumerate}
\end{proposition}

\noindent{\it Proof}.
We will just prove (\ref{difrulu}), the other relations being just as  simple.
If $g,h\in\cac_1 $, then
$$
\lc \der (gh) \rc_{st}
= g_th_t-g_sh_s
=g_s\lp h_t-h_s \rp +\lp  g_t-g_s\rp h_t\\
=g_s \lp \der h \rp_{st}+ \lp \der g \rp_{st} h_t,
$$
which proves our claim.

\fin
\vspace{0.3cm}

The iterated integrals of smooth functions on $\ott$ are obviously
particular cases of elements of $\cac$, which will be of interest for
us. Let us recall  some basic  rules for these objects:
consider $f\in\cac_1^\infty\cl^{d,l}$ and $g\in\cac_1^\infty(\R^d)$,
where $\cac_1^\infty $ denotes the set of
smooth functions on $\ott$. Then the integral $\int f \, dg$,
which will be denoted by
$\cj(f\, dg)$, can be considered as an element of
$\cac_2^\infty(\R^l)$. Namely, for $s,t\in\ott$ we set
$$
\cj_{st}(f\,  dg)
=
\left(\int  f dg \right)_{st} = \int_s^t   f_u dg_u.
$$
The multiple integrals can also be defined in the following way:
given a smooth element $h \in \cac_2^\infty\cl^{d,l}$ and $s,t\in\ott$, we set
$$
\cj_{st}(h\, dg )\equiv
\left(\int h dg  \right)_{st} = \int_s^t  h_{su} dg_u .
$$
In particular,  for
$f^1\in\cac_1^\infty(\R^{d_1})$, $f^2\in\cac_1^\infty\cl^{d_1,d_2}$
and $f^3\in\cac_1^\infty\cl^{d_2,d_3}$ the double integral
$\cj_{st}( f^3\, df^2 df^1)$ is defined  as
$$
\cj_{st}( f^3\, df^2df^1)
=\lp \int f^3\, df^2 df^1  \rp_{st}
= \int_s^t \cj_{su}\lp f^3\, df^2  \rp \, df_u^1.
$$
Now suppose that the $n$th order iterated integral of
$f^{n+1}df^n\cdots df^2$, which is  denoted by
$\cj(f^{n+1}df^n\cdots df^2)$, has been defined for
$f^j\in\cac_1^\infty\cl^{d_{j-1},d_j}$.
Then, if $f^1\in\cac_1^\infty(\R^{d_1})$, we set
\begin{equation}\label{multintg}
\cj_{st}(f^{n+1}df^n \cdots df^2 df^1)
=
\int_s^t  \cj_{su}\lp f^{n+1}df^n\cdots df^2\rp  \, df_u^{1},
\end{equation}
which recursively defines the iterated integrals of smooth functions.
Observe that a $n$th order integral $\cj(df^n\cdots df^2 df^1)$ could be
defined along the same lines.

\medskip

The following relations between multiple integrals and the operator $\der$ will also be useful:
\begin{proposition}\label{dissec}
Let $f\in\cac_1^\infty\cl^{d,l}$ and $g\in\cac_1^\infty(\R^d)$.
Then it holds that
$$
\der g = \cj( dg), \qquad
\der\lp \cj(f dg)\rp = 0, \qquad
\der\lp \cj (df dg)\rp = - (\der f) (\der g) = -\cj(df) \cj(dg),
$$
and
$$
 \der \lp \cj( df^n \cdots df^1)\rp  =
- \sum_{i=1}^{n-1}
\cj\lp df^n \cdots df^{i+1}\rp \cj\lp df^{i}\cdots df^1\rp.
$$
\end{proposition}

\noindent{\it Proof}.
Here  the proof is elementary again. We will just show the third of
the relations.  For $s,t\in\ott$ we have
$$
\cj_{st} (dg df)
= \int_s^t (f_u-f_s) dg_u
= \int_s^t  f_u dg_u - K_{st},
$$
with $K_{st}=f_s(g_t-g_s)$. The first term of the right hand side is easily seen to be in $\cz\cac_2$. Thus
$$
\der\lp \cj (dg df)\rp_{sut} =-
\lp\der K\rp_{sut} = -[f_u-f_s][g_t-g_u],
$$
which gives the announced result.
\fin

\subsection{Weakly controlled processes}\label{sec:wc-ps}

Recall that we have in mind to solve equations of the form
\begin{equation}\label{eq:rough-eds}
dy_t=\si(y_t) dx_t,\qquad y_0=a,
\end{equation}
where $t\in\ott$, $y$ is a $\R^l$-valued continuous process,
$\si:\R^l\to\cl^{d,l}$ is a $C_b^2$ function, i.e. twice continuously differentiable and bounded together with its derivatives, $x$ is a $\R^d$-valued
path and $a\in\R^l$ is a fixed initial
condition. As usual in rough path type considerations, we will have to assume
a  priori the following hypothesis  in order
to handle equations like (\ref{eq:rough-eds}):
\begin{hypothesis}\label{hyp:x}
The path $x$ is $\R^d$-valued
$\ga$-H\"older with $\ga>1/3$ and  admits a L\'evy area,
that is a process $\xd=\cj(dx dx)\in\cac_2^{2\ga}\cl^{d,d}$
satisfying
$$
\der\xd=\der x\otimes \der x,
\quad\mbox{i.\!\! e.}\quad
\lc (\der\xd)_{sut} \rc(i,j)
=
[\der x^{i}]_{su} [\der x^{j}]_{ut},
\quad s,u,t\in\ott, \, i,j\in\{1,\ldots,d  \}.
$$
\end{hypothesis}
\noindent
The solution to (\ref{eq:rough-eds}) will then be
expressed as a continuous function
of the input $a,\si,x$ and $\xd$.

\vspace{0.3cm}

Let us now be  more specific about the global strategy, we will use
to solve equation (\ref{eq:rough-eds}). First of all,  simple heuristic
considerations show that, if the equation admits a solution, it
should be a weakly controlled path, i.e.,  a process of the
following form:
\begin{definition}\label{def39}
Let $z$ be a process in $\cac_1^\ka(\R^k)$ with $\ka\le\ga$
and $2\ka+\ga>1$.
We say that $z$ is a weakly controlled path based on $x$, if
$z_0=a$, which is a given initial condition in $\R^k$,
and $\der z\in\cac_2^\ka(\R^k)$ can be decomposed into
\begin{equation}\label{weak:dcp}
\der z=\zeta \der x+ r,
\quad\mbox{i.\!\! e.}\quad
(\der z)_{st}=\zeta_s (\der x)_{st} + r_{st},
\quad s,t\in\ott,
\end{equation}
with $\zeta\in\cac_1^\ka\cl^{d,k}$ and $r$ is a regular part
such that $r\in\cac_2^{2\ka}(\R^k)$. The space of weakly controlled
paths will be denoted by $\cq_{\ka,a}(\R^k)$, and a process
$z\in\cq_{\ka,a}(\R^k)$ can be considered in fact as a couple
$(z,\zeta)$. The natural semi-norm on $\cq_{\ka,a}(\R^k)$ is given
by
$$
\cn[z;\cq_{\ka,a}(\R^k)]=
\cn[z;\cac_1^{\ka}(\R^k)]
+ \cn[\zeta;\cac_1^{\infty}\cl^{d,k}]
+ \cn[\zeta;\cac_1^{\ka}\cl^{d,k}]
+\cn[r;\cac_2^{2\ka}(\R^k)]
$$
with $\cn[g;\cac_1^{\ka}]$ defined by
(\ref{def:hnorm-c1}) and
$\cn[\zeta;\cac_1^{\infty}(V)]=\sup_{0\le s\le T}|\zeta_s|_V$.
\end{definition}

Note that it is always possible to find $\kappa \leq \gamma$ with $2 \kappa + \gamma >1$, since $\ga>1/3$. 
With this definition at  hand, we will try to solve equation
(\ref{eq:rough-eds}) in the following way:
\begin{enumerate}
\item
Study the stability of $\cq_{\ka,a}(\R^k)$ under a smooth map
$\vp:\R^k\to\R^n$.
\item
Define rigorously the integral $\int z_u dx_u=\cj(z dx)$
for a weakly controlled path $z$ and compute its decomposition
(\ref{weak:dcp}).
\item
Solve equation (\ref{eq:rough-eds}) in the space $\cq_{\ka,a}(\R^k)$
by a fixed point argument.
\end{enumerate}
In this section, we will concentrate on the first two points of this
program.

\vspace{0.3cm}

Let us first see, how smooth functions act on weakly controlled paths:
\begin{proposition}\label{cp:weak-phi}
Let $z\in\cq_{\ka,a}(\R^k)$ with decomposition (\ref{weak:dcp}),
$\vp \in   C_{b}^{2}(\R^k;\R^n)$ and set $\hz=\vp(z)$, $\ha=\vp(a)$.
Then $\hz\in\cq_{\ka,\ha}(\R^n)$, and it can be decomposed into
$$
\der \hz= \hat\zeta \der x +\hr,
$$
with
$$
\hat\zeta= \nabla\vp(z)\zeta
\quad\mbox{ and }\quad
\hr= \nabla\vp(z) r + \lc \der(\vp(z))-\nabla\vp(z)\der z \rc.
$$
Furthermore,
\begin{equation}\label{bnd:phi}
\cn[\hz;\cq_{\ka,\ha}(\R^n)]\le
c_{\vp,T}\lp 1+\cn^2[z;\cq_{\ka,a}(\R^n)]  \rp.
\end{equation}
\end{proposition}
\noindent{\it Proof}.
The algebraic part of the assertion is quite straightforward. Just write
\begin{eqnarray*}
(\der\hz)_{st}&=&\vp(z_t)-\vp(z_s)=
\nabla\vp(z_s)(\der z)_{st}+ \vp(z_t)-\vp(z_s)-\nabla\vp(z_s)(\der z)_{st}\\
&=&\nabla\vp(z_s)\zeta_s(\der x)_{st}
+\nabla\vp(z_s) r_{st}+ \vp(z_t)-\vp(z_s)-\nabla\vp(z_s)(\der z)_{st}\\
&=& \hat\zeta_s (\der x)_{st} + \hr_{st},
\end{eqnarray*}
which is the desired decomposition.

\vspace{0.3cm}

In order to give an estimate for $\cn[\hz;\cq_{\ka,\ha}(\R^n)]$, one has of course
to establish bounds for $\cn[\hz;\cac_1^{\ka}(\R^n)]$, $\cn[\hat\zeta;$
$\cac_1^{\ka}\cl^{d,n}]$, $\cn[\hat\zeta;\cac_1^{\infty}\cl^{d,n}]$
and $\cn[\hr;\cac_2^{2\ka}(\R^n)]$. Let us focus on the last of these
estimates, the other ones are quite similar. First notice that
$\hr=\hr^1+\hr^2$ with
\begin{equation}\label{def:hr12}
\hr^1_{st}=\nabla\vp(z_s) r_{st}
\quad\mbox{ and }\quad
\hr^2_{st}=\vp(z_t)-\vp(z_s)-\nabla\vp(z_s)(\der z)_{st}.
\end{equation}
Now, since $\nabla\vp$ is a bounded $\cl^{k,n}$-valued function, we have
\begin{equation}\label{est:hr1}
\cn[\hr^1;\cac_2^{2\ka}(\R^n)]
\le
\|\nabla\vp\|_\infty \cn[r;\cac_2^{2\ka}(\R^k)].
\end{equation}
Moreover,
\begin{equation*}
|\hr^2_{st}|\le \frac12 \|\nabla^2\vp\|_\infty |(\der z)_{st}|^2
\le c_{\vp} \cn^{2}[z;\cac_1^\ka(\R^k)]|t-s|^{2\ka},
\end{equation*}
which yields
\begin{equation}\label{est:hr2}
\cn[\hr^2;\cac_2^{2\ka}(\R^n)]\le
c_{\vp} \cn^2[r;\cac_2^{2\ka}(\R^k)],
\end{equation}
and thus
we obtain
$$
\cn[\hr;\cac_2^{2\ka}(\R^n)]\le
c_{\vp} \lp 1+\cn^2[r;\cac_2^{2\ka}(\R^k)]\rp,
$$
which ends the proof.
\fin

Let us now turn  to the integration of weakly controlled paths, which is
summarized in the following proposition. Notice that below we will use
two additional notations. We will set $M^*$ for the transposition of a
matrix $M$ and denote by $M\cdot N$ the inner product of two vectors
or two matrices.
\begin{proposition}\label{intg:mdx}
For a given $\ga>1/3$ and $\ka<\ga$,
let $x$ be a process satisfying Hypothesis \ref{hyp:x}. Furthermore,  let
$m\in\cq_{\ka,b}(\cl^{d,1})$ with decomposition $m_0=b\in\cl^{d,1}$
and
\begin{equation}\label{dcp:m}
(\der m)_{st}=\lc\mu_s (\der x)_{st}\rc^* +r_{st},
\quad\mbox{ where }\quad
\mu\in\cac_1^\ka\cl^{d, d}, \, r\in\cac_2^{2\ka}\cl^{d,1}.
\end{equation}
Define $z$ by $z_0=a\in\R$ and
\begin{equation}\label{dcp:mdx}
\der z=
m\,\der x+ \mu \cdot \xd \mathbf{-}\laa(r\, \der x+\der \mu\cdot \xd).
\end{equation}
Finally, set
 \begin{equation} \cj(m\, dx) = \der z. \end{equation}
Then:
\begin{enumerate}
\item
$z$ is well-defined as an element of $\cq_{\ka,a}(\R)$.
\item
The semi-norm of $z$ in $\cq_{\ka,a}(\R)$ can be estimated as
\begin{equation}\label{bnd:norm-imdx}
\cn[z;\cq_{\ka,a}(\R)]\le
c_{x}
\lp 1 + T^{\ga-\ka}\cn[m;\cq_{\ka,b}(\cl^{d,1})]\rp,
\end{equation}
for a positive constant $c_{x}$ depending only on $x$ and $\xd$.
The constant $c_{x}$ can be bounded as follows:
$$
c_x\le c
\lc
|x|_{\ga}+|\xd|_{2\ga}
\rc,
\quad\mbox{ for a universal constant }c.
$$
Moreover,
we have
\begin{equation}\label{bnd:inc-z}
\| \der z \|_{\ka}\le c_{x} T^{\ga-\ka}\cn[m;\cq_{\ka,b}(\cl^{d,1})].
\end{equation}
\item
It holds
\begin{equation}\label{rsums:imdx}
\cj_{st}(m\, dx)
=\lim_{|\Pi_{st}|\to 0}\sum_{i=0}^n
\lc m_{t_{i}}(\der x)_{t_{i}, t_{i+1}}
+ \mu_{t_{i}} \cdot\xd_{t_{i}, t_{i+1}} \rc
\end{equation} for any $0\le s<t\le T$,
where the limit is taken over all partitions
$\Pi_{st} = \{s=t_0,\dots,t_n=t\}$
of $[s,t]$, as the mesh of the partition goes to zero.
\end{enumerate}
\end{proposition}

\vspace{0.3cm}

Before going into the technical details of the proof,  let us see
how to recover (\ref{dcp:mdx}) in the smooth case, in order to
justify our  definition of the integral. (Notice however that
(\ref{rsums:imdx}) corresponds to the usual definition in the rough paths theory
\cite{LyonsBook}, which gives another kind of justification.)
\\ Let us assume for the moment
that $x$ is a smooth function and that $m\in\cac_1^\infty(\cl^{d,1})$
admits the decomposition (\ref{dcp:m}) with
$\mu\in\cac_1^\infty\cl^{d,d}$ and $r\in\cac_2^{\infty}\cl^{d,1}$.
Then $\cj(m \, dx)$ is well-defined, and we have
$$
\ist m_u dx_u = m_s[x_t-x_s] + \ist[m_u-m_s] dx_u
$$  for $s\le t$,
respectively
$$
\cj(m \,dx)= m\,\der x + \cj(\der m \, dx).
$$
Let us now plug the  decomposition (\ref{dcp:m}) into this expression,
which yields
\begin{eqnarray}\label{eq:imdx}
\cj(m \,dx)&=& m\,\der x + \cj((\mu \der x)^*\, d x) + \cj(r \,dx)\nonumber\\
&=& m\,\der x +\mu\cdot \xd + \cj(r \,dx).
\end{eqnarray}
For the sake of clarity, let us give some details about the identity
$\cj((\mu \der x)^*\, d x)=\mu\cdot \xd$. Indeed, according to our definitions
in Section \ref{cpss} we have
$$
\cj_{st}((\mu \der x)^*\, d x)
= \ist  \lp [\der x]_{su}^* \mu_s \rp dx_u
= \ist \mu_s \cdot \lc (\der x)_{su} \otimes dx_u  \rc
= \mu_s \cdot \ist   [\der x]_{su} \otimes dx_u
=\mu_s \cdot \bx_{st}^{\bf 2}
$$ for $s\le t$,
which proves the announced identity.
Notice also that the terms
$m\,\der x$ and $\mu\, \xd$  in (\ref{eq:imdx})  are well-defined as soon as $x$ and $\xd$
are defined themselves. In order to push forward our analysis to the
rough case, it remains to handle the term $\cj(r \,dx)$.
Thanks to (\ref{eq:imdx}) we can write
$$
\cj(r \,dx)= \cj(m \,dx)- m\,\der x -\mu\cdot \xd,
$$
and let us analyze this relation by applying $\der$ to both sides.
Using the second part of Proposition \ref{difrul} and the whole Proposition \ref{dissec}
yields
\begin{eqnarray}\label{eq:di-rdx}
\der\lc\cj(r \,dx)\rc&=&
- \der\lc m\,\der x\rc -\der\lc \mu\cdot \xd\rc\nonumber\\
&=& -\der m\, \der x -\der\mu\cdot\xd
+\mu\cdot(\der x\otimes\der x)\nonumber\\
&=&-\lc (\mu \,\der x)^* + r \rc \der x
-\der\mu\cdot\xd+\mu\cdot(\der x\otimes\der x)\nonumber\\
&=& -\der\mu\cdot\xd - r\, \der x.
\end{eqnarray}
Assuming now that $\der\mu\cdot\xd$ and $r\, \der x$ are both elements
of $\cac_2^\mu$ with $\mu>1$, $\der\mu\cdot\xd + r\, \der x$ becomes
an element of $\dom(\laa)$, and inserting (\ref{eq:di-rdx}) into
(\ref{eq:imdx}) we obtain
$$
\der z=\cj(m\, dx)\equiv
m\,\der x+ \mu \cdot\xd -\laa(r\, \der x+\der \mu\cdot \xd),
$$
which is the expression (\ref{dcp:mdx}) of our Proposition \ref{intg:mdx}.
Thus (\ref{dcp:mdx}) is a natural expression
for $\cj(m\,  dx)$.

\vspace{0.5cm}

\noindent
{\it Proof of Proposition \ref{intg:mdx}}.
We will decompose this proof in two steps.

\vspace{0.3cm}

\noindent
{\it Step 1:}
Recalling the assumption $2\ka+\ga>1$,
let us analyze the three terms in the right hand side of (\ref{dcp:mdx})
and show that they define an element of $\cq_{\ka,a}$ such that
$\der z =\zeta\der x +r$ with
$$
\zeta=m
\quad\mbox{ and }\quad
\hat r= \mu\cdot \xd- \laa\lp r\,\der x+\der \mu\cdot \xd \rp.
$$
Indeed, on the one hand $m\in\cac_1^\ka\cl^{d,1}$ and thus $\ze=m$ is of
the desired form for an element of $\cq_{\ka,a}$. On the other hand,
if $m\in\cq_{\ka,b}$, $\mu$ is assumed to be bounded and since
$\xd\in\cac_2^{2\ka}\cl^{d,d}$ we get that
$\mu\cdot\xd\in\cac_2^{2\ga}(\R)$.
Along the same lines we can prove that $r\,\der x\in \cac_3^{2\ka+\ga}(\R)$
and $\der \mu\cdot\xd\in \cac_3^{\ka+2\ga}(\R)$. Since
$\ka+2\ga\ge 2\ka+\ga>1$, we obtain  that
$r\,\der x+\der \mu\cdot\xd\in\dom(\laa)$
and
$$
\laa\lp r\,\der x+\der \mu\,\xd \rp\in \cac_2^{2\ka+\ga}(\R).
$$
Thus we have proved that
$$
r= \mu\cdot \xd- \laa\lp r\,\der x+\der \mu\cdot \xd \rp \in
\cac_2^{2\ka}(\R)
$$
and hence that $z\in\cq_{\ka,a}(\R)$. The estimates (\ref{bnd:norm-imdx})
and (\ref{bnd:inc-z})
are now obtained using  to the same kind of considerations and are left
to the reader for the sake of conciseness.

\vspace{0.3cm}

\noindent
{\it Step 2:}
The same kind of computations as those leading to (\ref{eq:di-rdx})
also show that
$$
\der\lp m\, \der x+\mu \cdot\xd  \rp
=
-\lc r\, \der x+ \der\mu\cdot \xd \rc.
$$
Hence equation (\ref{dcp:mdx}) can also be read as
$$
\cj(m\, dx)=\lc \id-\laa\der \rc (m\, \der x+\mu \cdot\xd ),
$$
and a direct application of Corollary \ref{cor:integration} yields
(\ref{rsums:imdx}), which ends our proof.

\fin

\vspace{0.3cm}

Notice that the previous proposition has a straightforward multidimensional
extension, which we state for further use:
\begin{corollary}
Let $x$ be a process satisfying Hypothesis \ref{hyp:x} and let
$m\in\cq_{\ka,b}(\cl^{d,k})$ with decomposition $m_0=b\in\cl^{d,k}$
and
\begin{equation}\label{dcp:m-mult}
(\der m^i)_{st}=\lc\mu_s^i (\der x)_{st}\rc^* +r_{st}^i
\quad\mbox{ with }\quad
\mu^i\in\cac_1^\ka\cl^{d, d}, \, r^i\in\cac_2^{2\ka}\cl^{d,k}, \quad  i=1, \ldots,  k,
\end{equation}
where we have considered $m$ as a $\R^{k\times d}$-valued path
and have set $m^i=m(i,\cdot)$.
Define $z$ by $z_0=a\in\R^k$ and
\begin{equation}\label{dcp:mdx-mult}
\der z^i=\cj(m^i\, dx)\equiv
m^i\,\der x+ \mu^i \cdot \xd +\laa(r^i\, \der x+\der \mu^i\cdot \xd).
\end{equation}
Then the conclusions of Proposition \ref{intg:mdx} still hold  in
this context.
\end{corollary}

Notice also that our extended pathwise integral has a nice
continuity property with respect to the driving path $x$.
See also \cite[p. 14]{Gu}.
\begin{proposition}\label{cvgce:wc-ps}
Let $x$ be a function satisfying Hypothesis \ref{hyp:x} and
assume that there exists a sequence $\{x^n;\, n\ge 1  \}$
of piecewise $C^1$-functions from $\ott$ to $\R^d$ such that
\begin{equation}\label{hyp:lim-levy}
\lim_{n\to\infty} \cn[x^n-x;\cac_1^{\ga}(\R^d)]=0,
\quad\mbox{ and }\quad
\lim_{n\to\infty} \cn[\bx^{{\bf 2},n}-\xd;\cac_2^{2\ga}\cl^{d,d}]=0.
\end{equation}
For $n\ge 1$, define $z^n\in\cac_1^\ka(\R^k)$ in the
following way: set $z_0=b\in\R^k$ and
$$
\der z^n =\zeta^n \der x^n + r^n,
$$
where $\zeta^n\in\cac_1^\ka\cl^{d,k}$ and $r^n\in\cac_2^{2\ka}(\R^k)$.
Moreover, let $z$ be a weakly controlled process with  decomposition
(\ref{weak:dcp}) and assume that
$$
\lim_{n\to\infty}
\cn[z^n-z;\cac_1^{\ka}(\R^k)]+
\cn[\zeta^n-\zeta;\cac_1^{\infty}\cl^{d,k}]+
\cn[\zeta^n-\zeta;\cac_1^{\ka}\cl^{d,k}]+
\cn[r^n-r;\cac_2^{2\ka}(\R^k)]
=0.
$$
Eventually, let $\vp:\R^k\to\cl^{d,m}$ be a $C_b^2$-function. Then
$$
\lim_{n\to\infty}
\cn\lc \cj(\vp(z^n) dx^n)- \cj(\vp(z) dx);\, \cac_2^{\ka}(\R^m)\rc
=0.
$$
\end{proposition}

\subsection{Stochastic calculus with respect to a rough path}
\label{sec:stoch-rough}

In this section, we will apply the previous considerations to
two of the usual main aims in the theory of  stochastic calculus:
to study differential equations driven by a rough
signal and to establish a  change of variable  formula.

\subsubsection{Rough differential equations}

Recall that we wish to solve equations of the form (\ref{eq:rough-eds}).
In our algebraic setting, we will rephrase this as follows: we will
say that $y$ is a solution to (\ref{eq:rough-eds}), if $y_0=a$,
$y\in\cq_{\ka,a}(\R^l)$ and for any $0\le s\le t\le T$ we have
\begin{equation}\label{eds:alg-form}
(\der y)_{st}=\cj_{st}(\si(y)\, dx),
\end{equation}
where the integral $\cj(\si(y)\, dx)$ has to be understood in the sense
of Proposition \ref{intg:mdx}. Our existence and uniqueness result
reads as follows:
\begin{theorem}\label{thm:ex-uniq}
Let $x$ be a process satisfying Hypothesis \ref{hyp:x} and
$\si:\R^l\to\cl^{d,l}$ be a $C^2$ function, which is bounded
together with its derivatives. Then
\begin{enumerate}
\item
Equation (\ref{eds:alg-form}) admits a unique solution $y$ in
$\cq_{\ka,a}(\R^l)$ for any $\ka<\ga$ such that $2\ka+\ga>1$.
\item
The mapping $(a,x,\xd)\mapsto y$ is continuous from
$\R^l\times\cac_1^{\ga}(\R^d)\times\cac_2^{2\ga}(\R^{d\times d})$
to $\cq_{\ka,a}(\R^l)$.
\end{enumerate}
\end{theorem}

\noindent{\it Proof}.
We will identify the solution on a small interval $[0,\tau]$ as
the fixed point of the map $\Gamma:\cq_{\ka,a}(\R^l)\to\cq_{\ka,a}(\R^l)$
defined by $\Gamma(z)=\hz$ with $\hz=a$ and $\der\hz=\cj(\si(z)\, dx)$.
The first step in this direction is to show that the ball
\begin{equation}\label{def:ball-m}
B_M=\lcl z; \,z_0=a, \,\cn[z;\cq_{\ka,a}(\R^l)]\le M \rcl
\end{equation}
is invariant under $\Gamma$ if $\tau$ is small enough and
$M$ is large enough. However, due to  Propositions
\ref{cp:weak-phi} and \ref{intg:mdx} and assuming $\tau \leq 1$
we have
\begin{equation}\label{bnd:ggaz}
\cn[\Gamma(z);\cq_{\ka,a}]\le
c_{\si,x}\lp 1+\tau^{\ga-\ka}\cn^2[z;\cq_{\ka,a}] \rp.
\end{equation}
Since the set ${\cal A}=\{u\in\R_+^*:\,\,c_{\si,x}( 1+\tau^{\ga-\ka}u^2)\le u\}$
is not empty as soon as $\tau$ is small enough
(see also the third point of the proof of Lemma \ref{ctrl-lm} below),
it is easily
shown that the ball $B_M$ defined at (\ref{def:ball-m}) is left
invariant by $\Gamma$ for $\tau$ small enough and $M$ in ${\cal A}$.

\vspace{0.3cm}

Now, since we are working in $B_M$, the fixed point argument for $\Gamma$
is a  standard argument and is left to the reader. This
leads to a unique solution to equation (\ref{eds:alg-form}) on a
small interval $[0,\tau]$. One is then able to obtain the unique solution on
an arbitrary interval $[0,k\tau]$ with $k\ge 1$ by patching solutions
on $[j\tau,(j+1)\tau]$. Notice here that an important point,
which allows us to use a constant step $\tau$, is the fact that
the estimate (\ref{bnd:ggaz}) does not depend on the initial condition
$a$, due to the fact that $\si$ is bounded together with its derivatives.

\fin

\begin{remark}{\rm
The case of an equation of the form
$$
d y_t = b(y_t) \, dt +\si(y_t)\, dx_t,
\quad \quad
y_0=a,
$$
which can be written equivalently
\begin{equation}\label{rough:drift}
\der y= \cj( b(y) \, dh) + \cj(\si(y)\, dx),  \quad \quad y_0=a,
\end{equation}
where $b:\R^k\to\R^k$ is a $C_b^1$-function and
$h:\ott\to\R$ is defined by $h_t=t$, can be solved easily
thanks to the previous results by taking one of the following
two observations into account:
\begin{enumerate}
\item
One can define the integral $\cj( b(y) \, dh)$ in the usual Riemann
sense for a weakly controlled path $y$. Then the fixed point argument of
Theorem \ref{thm:ex-uniq} can be extended trivially to the case of
an equation with drift.
\item
One can also define a path $\tilde x=(x,h^{(k)})$, where $h^{(k)}$
represents $k$ copies of the function $h$, and write equation
(\ref{rough:drift}) as $\der y =\cj(\tilde\si(y)\, d\tilde x)$,
for a new matrix-valued function $\tilde\si$. The existence and
uniqueness result follows then directly from Theorem \ref{thm:ex-uniq}
if $b$ is a $C_b^2$-function.
\end{enumerate}}
\end{remark}

\noindent

\begin{remark}{\rm
We have stressed here the fact that one could deal in a rather elementary way
with processes having a regularity $\ga>1/3$. However, if $\ga>1/2$,
our algebraic setting also applies and the results we have obtained
so far can be expressed in a simpler way:

\vspace{0.5cm}

Let $x$ be a $\R^d$-valued $\ga$-H\"older function with $\ga>1/2$,
and $m$ a function in $\cac_1^k\cl^{d,k}$, with $\ga+\ka>1$.
Define $z$ by $z_0=a\in\R^d$ and
\begin{equation}
\der z=
m\,\der x-\laa(\der m\, \der x),
\end{equation}
and set
 \begin{equation}
\cj(m\, dx) = \der z.
\end{equation}
Then:
\begin{enumerate}
\item
$z$ is well-defined as an element of $\cac_2^\ga(\R^k)$, and it holds:
\begin{equation}\label{est_rs}
\lln \cj_{st}(m\, dx)  \rrn
\le
\| m\|_{\infty} \|x\|_{\ga} |t-s|^{\ga}
+ c_{\ga,\ka} \| m\|_{\ka} \| x\|_{\ga} |t-s|^{\ga+\ka}.
\end{equation}
\item
The integral $\cj(m dx)$ coincides with the usual Young integral, and in
particular, it holds
\begin{equation}
\cj_{st}(m\, dx)
=\lim_{|\Pi_{st}|\to 0}\sum_{i=0}^n
m_{t_{i}}(\der x)_{t_{i}, t_{i+1}}
\end{equation} for any $0\le s<t\le T$,
where the limit is taken over all partitions
$\Pi_{st} = \{s=t_0,\dots,t_n=t\}$
of $[s,t]$, as the mesh of the partition goes to zero.
\item
Equation (\ref{eds:alg-form}) can be solved in the Young sense whenever
$\si$ is a $C_b^2$-function, and the solution $y$ satisfies
\begin{equation}\label{est_sol_rp}
\|y\|_{\la}\le c_{a,\si,T} \| x\|_{\ga},
\end{equation}
for all $\la<\ga$.
\end{enumerate}}
\end{remark}

\subsubsection{It\^o's type formula}

In the sequel of this paper, it will also be essential to have a
change of variable formula for weakly controlled process.
This will be achieved under the following additional assumption
on $\xd$, which will be shown to be valid in the fractional Brownian
motion case:
\begin{hypothesis}\label{hyp:sym-xd}
Let $\xd$ be the area process defined in Hypothesis \ref{hyp:x}
and denote by $\xds$ the symmetric part of $\xd$, i.e.
$\xds=\frac12(\xd+(\xd)^*)$. Then we assume that for $0\le s<t\le T$
we have
$$
\xsst= \frac12[\der x]_{st}\otimes [\der x]_{st}.
$$
\end{hypothesis}

\begin{remark}{\rm
It is worth noticing at this point that this assumption does not
involve any limit type property of the form (\ref{hyp:lim-levy})
for $\xd$. This will simplify the verification of the different hypothesis
for the fractional Brownian motion
with respect to \cite{CQ,LyonsBook}.}
\end{remark}

With these assumptions in mind, our change of variable formula  reads
as follows:
\begin{proposition}
Assume that $x$ satisfies Hypothesis \ref{hyp:x} and \ref{hyp:sym-xd}.
Let $m\in\cq_{\ka,b}(\cl^{d,k})$ be a process of the form
(\ref{dcp:m}) and let $z\in\cq_{\ka,a}(\R^k)$ be defined by $z_0=a$
and $\der z=\cj(m \,dx)$, which is given by (\ref{dcp:mdx-mult}). Let also
$f \in  C_b^2(\R^{k};\R)$. Then $f(z_0)=f(a)$ and
\begin{equation}\label{itoform}
\lc \der(f(z)) \rc_{st}= \cj_{st}\lp (\nabla^* f(z) m)\, dx \rp.
\end{equation}
\end{proposition}

\noindent{\it Proof}.
The strategy of our proof is quite straightforward. By using the
composition and integration rules for weakly controlled processes
we will compute the decompositions of $\der(f(z))$
and $\cj( (\nabla^* f(z) m)\, dx )$ respectively, and then show that
they coincide. Let us begin with the decomposition
of $\der(f(z))$. Recall that the decomposition (\ref{dcp:m-mult})
of $\der m$ can be written as
$$
\der m^i=\lc \mu^i\, \der x\rc^* + r^i,
\quad \quad   i=1, \ldots, k,
$$
and that we have $\der z^i=\cj(m^i\, dx)$. Thus, for $i\le k$
the decomposition of $z^i$ is given by
$$
\der z^i = m^i\, \der x + \mu^i\cdot \xd -
\laa\lp r^i\, \der x + \der\mu^i\cdot \xd \rp.
$$
In the sequel of the proof, we will use the following notation:
we write $\hr$ for any increment in $\cac_2^{\mu}$ with
$\mu\ge 2\ka+\ga>1$, whose
exact expression can change from line to line. In the same spirit, we will denote
by $\hrk$ any increment in $\cac_2^{2\ka}$ with regularity
at least $2\ka$. With these conventions in mind,
some elementary algebraic manipulations yield
\begin{align}\label{dcp:fz}
&\lc \der f(z)  \rc_{st}\nonumber\\
&=
\sum_{i=1}^{k}\partial_i f(z_s) \lc \der z^i  \rc_{st}
+\frac12 \sum_{i,j=1}^k
\partial_{ij}^{2}f(z_s) \lc \der z^i  \rc_{st} \lc \der z^j \rc_{st}
+ \hr_{st}\nonumber\\
&= \sum_{i=1}^{k}\partial_i f(z_s)
\lp m_s^i\, (\der x)_{st} + \mu_s^i\cdot \xdst \rp
+\frac12 \sum_{i,j=1}^k \partial_{ij}^{2}f(z_s)
\lc m_s^i (\der x)_{st} m_s^j (\der x)_{st}   \rc + \hr_{st} \\
&= \sum_{i=1}^{k}\partial_i f(z_s)
\lp m_s^i\, (\der x)_{st} + \mu_s^i\cdot \xdst \rp
+\frac12 \sum_{i,j=1}^k \partial_{ij}^{2}f(z_s)
\lc (m_s^i)^* m_s^j\rc\cdot \lc (\der x)_{st}\otimes (\der x)_{st} \rc
+ \hr_{st}.
\nonumber
\end{align}
which is the decomposition we were looking for $\delta f (z)$.

\vspace{0.3cm}

Let us compute now the decomposition of $\der[\nabla^*f(z) m]$.
We have due  to Proposition \ref{difrul} that
\begin{equation}\label{eq:nablaf-m}
\der[\nabla^*f(z) m]_{st}
=
\sum_{i=1}^{k} \der\lc\partial_i f(z_s)\rc_{st} m_t^i
+ \partial_i f(z_s) \lc \der m^i \rc_{st}.
\end{equation}
Recall also that, setting $m=[m^1,\ldots,m^k]^*$, $\der z$ can
be decomposed into
$\der z=m\der x+\hrk$. Thus, according to Proposition \ref{cp:weak-phi},
one gets
$$
\der g(z)=\sum_{j=1}^{k}\partial_j g(z) \, m^j\, \der x +\hrk,
$$
for any smooth function $g:\R^k\to\R$. Plugging this equality into
(\ref{eq:nablaf-m}), we obtain that
$$
\der[\nabla^*f(z) m]_{st}= A_{st}+B_{st}+\hrk_{st},
$$
where
$$
A_{st}= \sum_{i,j=1}^k \partial_{ij}^k f(z_s)
m_s^j (\der x)_{st} m_s^i,
\quad\mbox{ and }\quad
B_{st}= \sum_{i=1}^{k} \partial_i f(z_s) \lc \mu_s^i (\der x)_{st} \rc^*.
$$
Notice that in the definition of $A_{st}$, $m_t^i$ has been replaced by
$m_s^i$, since the difference between the two expressions is again
a remainder $\hrk$. Now, a little elementary linear algebra shows
that
$$
A_{st}= (\der x)_{st}^* \sum_{i,j=1}^k \partial_{ij}^k f(z_s)
(m_s^j)^*  m_s^i,
$$
and hence the decomposition of $\nabla^*f(z) m$ as a weakly
controlled process can be written as
\begin{equation}\label{dcp:nabla-fz}
\der[\nabla^*f(z) m]
= \lc \nu \der x \rc^* +\hrk,
\end{equation}
with
$$
\nu_s=\sum_{i,j=1}^k \partial_{ij}^k f(z_s)(m_s^i)^*  m_s^j
+  \sum_{i=1}^{k} \partial_i f(z_s) \mu_s^i.
$$

\vspace{0.3cm}

With the  expression (\ref{dcp:nabla-fz}) at hand, we are now ready to
compute $\cj(\nabla^*f(z) m \, dx)$. Indeed, using Proposition
\ref{intg:mdx} we get
\begin{align}\label{dcp:cj-nabla-f}
&\cj_{st}(\nabla^*f(z) m \, dx)
=\nabla^*f(z_s) \,m_s\, (\der x)_{st}
+\nu_s\cdot\xdst+\hr_{st}\nonumber\\
&=\sum_{i=1}^{k}\partial_i f(z_s)
\lp m_s^i\, (\der x)_{st} + \mu_s^i\cdot \xdst \rp
+ \sum_{i,j=1}^k \partial_{ij}^{2}f(z_s)
\lc (m_s^i)^* m_s^j\rc\cdot \xdst + \hr_{st}.
\end{align}

\vspace{0.3cm}

If we now put  the expressions (\ref{dcp:fz})
and (\ref{dcp:cj-nabla-f}) together, we end up with
\begin{multline}\label{diff:ito}
[\der f(z)]_{st}-\cj_{st}(\nabla^*f(z) m \, dx)\\
=\sum_{i,j=1}^k \partial_{ij}^{2}f(z_s)
\lc (m_s^i)^* m_s^j\rc\cdot \lc \xdst -
\frac12\lc (\der x)_{st}\otimes (\der x)_{st} \rc \rc +\hr_{st}.
\end{multline}
Let us show that this last expression  depends only on the
symmetric part $\xds$ of $\xd$. Indeed, it is easily checked that,
if $H$ is a symmetric matrix of $\R^{d,d}$, $\{M^i;i\le k  \}$
a family of elements of $\R^{1,d}$ and $X\in\R^{d,d}$, then
\begin{multline*}
\sum_{i,j=1}^k H(i,j)\lc (M^i)^* M^j\rc \cdot X\\
=\sum_{i,j=1}^k \sum_{\al,\beta=1}^d
H(i,j) M^i(\al) M^j(\beta) X(\al,\beta)
=\sum_{i,j=1}^k H(i,j)\lc (M^i)^* M^j\rc \cdot X^{\bf s},
\end{multline*}
where $X^{\bf s}$ denotes the symmetric part of $X$. Applying this
identity to $H=\mbox{Hess}(f(x_s))$, $M^i=m_s^i$ and $X=\xdst$,
equation (\ref{diff:ito}) becomes
\begin{multline*}
[\der f(z)]_{st}-\cj_{st}(\nabla^*f(z) m \, dx)\\
=\sum_{i,j=1}^k \partial_{ij}^{2}f(z_s)
\lc (m_s^i)^* m_s^j\rc\cdot \lc \xsst -
\frac12\lc (\der x)_{st}\otimes (\der x)_{st} \rc \rc +\hr_{st}
=\hr_{st},
\end{multline*}
thanks to Hypothesis \ref{hyp:sym-xd}. We thus have shown that
\begin{equation}\label{diff:ito2}
\der f(z)-\cj(\nabla^*f(z) m \, dx)= \hr,
\end{equation}
for an  increment $\hr\in\cac_2^{\mu}$ with $\mu>1$.
Now  we are in the position to prove easily that $\hr=0$. If
we apply $\delta$ to the expression above, we find that $\hr\in\ker\delta$.
Thanks to Lemma \ref{exd}, there exists  a function $g\in\cac_1$
such that $\hr=\delta g$. Moreover, $g$ inherits the regularity of
$\hr$, and hence $g\in\cac_1^\mu$ with $\mu>1$, which means that
$g$ is a constant function and that $\hr=\der g=0$. Putting  these
considerations and equation (\ref{diff:ito2}) together, we finally get
$$
\der f(z)-\cj(\nabla^*f(z) m \, dx)=0,
$$
which finishes our proof.

\fin

\subsection{Application to the fractional Brownian motion}

All the previous constructions rely on the specific assumptions
we have made on the process $x$. In this
section, we will show that
the results
given at Sections \ref{sec:wc-ps} and \ref{sec:stoch-rough}
 apply to the fractional Brownian motion.

The combination of the following
proposition  with the results for the general theory allows us  to prove Theorem \ref{sum-fbm}.
\begin{proposition}\label{prop:hyp-fbm}
Let $B$ be a $d$-dimensional fractional Brownian motion
and suppose $H>1/3$.
Then almost all sample paths of  $B$ satisfy
the Hypothesis \ref{hyp:x} and \ref{hyp:sym-xd}.
\end{proposition}

\noindent{\it Proof}.
Let us first check Hypothesis \ref{hyp:x}. It is a classical fact
that $B\in\cac_1^\ga$ for any $1/3<\ga<H$, when $B$ is a fractional Brownian
motion with $H>1/3$. As far as $\xd$ is concerned, a natural choice
 is
$$
\xdst=\ist dB_u \otimes \int_s^u dB_v,
\quad\mbox{i. e.}\quad
\xdst(i,j)=\ist dB_u^i  \int_s^u dB_v^j,
\quad i,j\in\{1,\ldots,d  \},
$$ for $0\le s < t \le T$,
where the stochastic integrals are understood in the Stratonovich sense.
Then it is a classical result that $\xd$ is well-defined for
$H>1/3$ (see, e.g., \cite{PT2} for $i\ne j$
 and \cite{CN} for $i=j$).
The substitution formula 
for Stratonovich integrals also easily
yields  that $\der\xd=\der x\otimes\der x$. Furthermore, by
stationarity (\ref{stationarity}) and the scaling property (\ref{scaling}) of the fractional Brownian motion, we have that
$$
{\rm E}\lc |\xdst(i,j)|^2  \rc
= |t-s|^{4H} {\rm E}\lc |\bx_{01}^{\bf 2}(i,j)|^2  \rc
\le c |t-s|^{4H}.
$$
>From this inequality and thanks to the fact that $\xd$ is a process in the
second chaos of the fractional Brownian motion $B$, on which all
$L^p$ norms are equivalent for $p>1$, we get that
\begin{equation}\label{ineq:norm-xd}
{\rm E}\lc |\xdst(i,j)|^p  \rc
\le c_p |t-s|^{2pH}.
\end{equation}
In order to conclude that $\xd\in\cac_2^{2\ga}\cl^{d,d}$ for any
$\ga<1/3$, let us recall the following inequality from \cite{Gu}:
let $g\in \cac_2(V)$ for a given Banach space $V$; then,
for any $\ka>0$ and $p\ge 1$ we have
\begin{equation}\label{lem:garsia}
\| g\|_{\ka}\le c \lp U_{\ka+2/p;p}(g) + \| \der g\|_{\ga}\rp
\quad\mbox{ with }\quad
U_{\ga;p}(g)=
\lp \int_0^T\int_0^T \frac{|g_{st}|^p}{|t-s|^{\ga p}} \rp^p.
\end{equation}
By plugging inequality (\ref{ineq:norm-xd}) into (\ref{lem:garsia})
and recalling that $\der\xd=\der x\otimes\der x$, we obtain
 that $\xd\in\cac_2^{2\ga}\cl^{d,d}$ for any $\ga<H$, which shows that
$B$ satisfies Hypothesis \ref{hyp:x}.

\vspace{0.3cm}

The proof of Hypothesis \ref{hyp:sym-xd} is now a consequence of the
It\^o-Stratonovich formula for the fractional Brownian motion (see, e.g., \cite{AMN}).

\fin

\begin{remark}
{\rm
Proposition \ref{cvgce:wc-ps} implies
that the theory of rough paths presented here and the classical one of Lyons' type coincide
for the fractional Brownian motion with
$H>1/3$. In particular, consider the multiple integrals
$$
\ist dB_{u_1}^{\alpha_1}\int_{s}^{u_1}dB_{u_2}^{\alpha_2}
\cdots \int_{s}^{\alpha_n}dB_{u_{n-1}}^{\alpha_{n-1}},
\quad\mbox{ for }\quad
(\alpha_1,\ldots,\alpha_{n})\in \{0,\ldots,d  \}^n,
$$
with the convention $B_t^0=t$. Then these multiple integrals, which are constructed
by means of Proposition \ref{intg:mdx}, coincide with the usual Stratonovich
integral with respect to the fractional Brownian motion, see, e.g.,
\cite{AMN,BC}.
}
\end{remark}

\section{Malliavin calculus with respect to fBm}
In this section, we assume that the Hurst index of $B$ verifies $H>1/2$. Let us give a few facts about the Gaussian structure of fractional Brownian motion
and its Malliavin derivative process, following  Chapter 1.2 in \cite{nual} and Section 2 in \cite{nualsaus}.
Let $\EE$ be the set of step-functions on $[0,T]$ with values in $\R^{d}$.  Consider the Hilbert space $\HH$ defined as the closure of $\EE$ with respect to the scalar product induced by
$${ \left \lpa (\1_{[0,t_{1}]}, \ldots, \1_{[0,t_{d}]}), (\1_{[0,s_{1}]}, \ldots, \1_{[0,s_{d}]}) \right\rpa}_{\HH}\; =\; \sum_{i=1}^{d}R_H(t_{i},s_{i}), \quad s_{i},t_{i} \in [0,T], \,\, i=1, \ldots, d,$$
where $R_{H}(t,s)$ is given by (\ref{rhts}).
The scalar product between two elements $\phi, \psi \in \EE$ is given by
\begin{equation}{\label{sca}} \langle \varphi, \psi \rangle_{\HH} =\gamma_{H}  \sum_{i=1}^{d} \int_{0}^{1} \int_{0}^{1} \varphi^{i}(r) \psi^{i}(u) |r-u|^{2H-2} \, dr \, du \end{equation}
with $\gamma_{H}=H(2H-1)$. The space $\HH$ contains $L^{\frac{1}{H}}([0,T]; \R^{d})$ but its elements can be distributions, see, e.g., \cite{PT}. Formula (\ref{sca}) holds also for $\varphi, \psi \in L^{\frac{1}{H}}([0,T]; \R^{d})$. The mapping
$$ (\1_{[0,t_{1}]}, \ldots, \1_{[0,t_{d}]}) \mapsto \sum_{i=1}^{d}B_{t_{i}}^{i} $$can be extended to an isometry between $\HH$ and the Gaussian space $H_{1}(B)$ associated with $B=(B^{1}, \ldots , B^{d})$. We denote this isometry by $\varphi \mapsto B(\varphi)$.
Let $\mathcal{S}$ be the set of smooth cylindrical random variables of the form
$$F = f(B(\varphi_{1}), \ldots, B(\varphi_{k})), \qquad \varphi_{i} \in \HH, \quad i=1, \ldots, k,$$ where $f\in C^{\infty}(\R^{k},\R)$ is bounded with bounded derivatives. The derivative operator $D$ of a smooth cylindrical random variable of the above form is defined as the $\HH$-valued random variable
$$ DF= \sum_{i=1}^{k} \frac{\partial f}{\partial x_{i}} (B(\varphi_{1}), \ldots, B(\varphi_{n})) \varphi_{i}.$$
This operator is closable from $L^{p}(\Omega)$ into $L^{p}(\Omega; \HH)$.
 As usual, $\sk^{1,2}$ denotes
the closure of the set of smooth random variables with respect to the norm
$$\| F\|_{1,2}^2 \; = \; { \rm E} |F|^2  + { \rm E} \| D F\|_{\HH}^2 .$$
In particular, if $D^{i}F$ designates the Malliavin derivative
of  a functional $F\in\sk^{1,2}$ with respect to $B^i$,
we have  $D^{i} B^{j}_{t} = \delta_{i,j}\1_{[0,t]}$ for $i,j=1, \ldots, d$.

The divergence operator $\delta$ is the adjoint of the derivative operator. If a random variable $u \in L^{2}(\Omega; \HH)$
belongs to $\textrm{dom}(\delta)$, the domain of the divergence operator, then  $\delta(u)$ is defined by the duality relationship
\begin{equation}\label{eq:duality}
{ \rm E} (F \delta(u))= { \rm E} \langle D F, u \rangle_{\HH}, 
\end{equation}
for every $F \in \sk^{1,2}$.
Moreover, if $u \in \textrm{dom}(\delta)$ and $F \in \sk^{1,2}$ such that $Fu \in L^{2}(\Omega; \HH)$, then we have the following integration by parts formula
\begin{equation}\label{int_part}    \delta(Fu)=F\delta(u)- \langle DF, u \rangle_{\HH}. \end{equation}

The following proposition is well known.
For part (a) see, e.g., \cite{NR}
and for  part (b) and (c), see   Proposition 19 in \cite{nualsaus} and Theorem 3.1 in \cite{HN}.

\begin{proposition}{\label{prop_X}}  Let $b^{i}$ and $\sigma^{i,j}$, $i=1,  \ldots, n$, $j=1, \ldots, d$ be twice continuously differentiable with bounded derivatives. \\
(a) Then equation (\ref{eqintro}) has a unique solution $X=(X^{1}, \ldots, X^{n})$ in the Young sense in the class of all processes having $\alpha$-H\"older continuous sample paths with $1-H<\alpha<H$.  \\
(b) It holds
$$  \max_{i=1, \ldots, n} \, {\rm E}    \sup_{0\leq t \leq T} |X^{i}_{t}|^{p} \,    < \infty    $$
for all $p \geq 1$. \\
(c) Moreover, we have
$X^{i}_{t} \in \sk^{1,2}(\HH)$ for all $t\in [0,T]$, $i=1, \ldots, n$.  
The Malliavin derivative satisfies almost surely:
\begin{align*}
& D^{j}_{s}X_{t}^{i}= \sigma^{i,j}(X_{s})+ \sum_{k=1}^{n} \int_{s}^{t} b_{x_{k}}^{i}(X_{u})D_{s}^{j}X^{k}_{u} \, du + \sum_{k=1}^{n} \sum_{j=1}^{d} \int_{s}^{t} \sigma_{x_{k}}^{i,j}(X_{u}) D_{s}^{j}X^{k}_{u} \, dB^{j}_{u},
\quad s\le t,   \\ & D^{j}_{s}X_{t}^{i} = 0, \quad s>t,
\end{align*} for   $j=1, \ldots, d$,
 where $D^{j}_{s}X_{t}^{i}$ is the $j$-th component of $D_{s}X^{i}_{t}.$
Furthermore
\begin{eqnarray*}
  \max_{j=1, \ldots, d} \, \max_{i=1, \ldots, n} \, \sup_{0 \leq s \leq t \leq T} {\rm \E} |D^{j}_{s}X_{t}^{i} |^{p} < \infty.
\end{eqnarray*}
\end{proposition}

\section{Proof of Theorem \ref{thmptfa}}
In the present section we will prove Theorem \ref{thmptfa}. We
separate the proof in two parts: in the first one, we will show
how to use trees for the parametrization of  the expansion; while
in the second one we explain how to control the remainder term,
which appears when we expand $P_tf(a)$ with respect to $t$.

\subsection{Rooted trees approach}
In this section, we assume that the Hurst index of the fBm $B$ verifies $H>1/3$.
The first step in the proof of the algebraic part of Theorem \ref{thmptfa} is
the following result.
\begin{theorem} \label{St-tree-expansion-exact-sol:Wm}
    Let $(X_t^a)_{t \in [0,T]}$ be the solution of SDE~(\ref{eqintro})
    with initial value $X_{0}^a = a \in \mathbb{R}^n$.
    Then for $m \in \mathbb{N}_0$ and $f \in C_b^{m+1}(\mathbb{R}^n; \mathbb{R}),  b, \sigma^j,  \in C_b^{m+1}(\mathbb{R}^n; \mathbb{R}^n)$,
    $1 \leq j \leq d$, we get for $t \in [0,T]$ the expansion
    \begin{equation} \label{Strato-tree-expansion-exact-sol-formula1:Wm}
         f(X_t^a) = \sum_{\substack{\textbf{t} \in LTS \\
        l(\textbf{t})-1 \leq m}}
        \sum_{j_1, \ldots, j_{s(\textbf{t})}=1}^d
        F(\textbf{t})(a) \,\, \mathcal{I}_{\textbf{t}}(1)_{0,t}
        + \mathcal{R}_m(0,t)
    \end{equation}
    with a truncation term
    \begin{equation}
    \label{St-tree-expansion-exact-sol-truncation:Wm}
        \mathcal{R}_m(0,t) = \sum_{\substack{\textbf{t} \in LTS \\
        l(\textbf{t})-1 = m+1}}
        \sum_{j_1, \ldots, j_{s(\textbf{t})}=1}^d
        \mathcal{I}_{\textbf{t}}(F(\textbf{t})(X_s^a))_{0,t}
    \end{equation}
\end{theorem}
\quad \\
\noindent{\it Proof}. The proof is very similar to the proof of
Theorem~4.2 in~\cite{Roessler}.
Recall we defined the differential operators $\mathcal{D}^0$ and
$\mathcal{D}^j$ as
\begin{equation}
    \mathcal{D}^0 = \sum_{k=1}^n b^k \, \frac{\partial}{\partial x^k} \qquad \qquad
    \text{ and } \qquad \qquad \mathcal{D}^j = \sum_{k=1}^n \sigma^{k,j} \,
    \frac{\partial}{\partial x^k}
\end{equation}
for $j=1, \ldots, d$.
Moreover, set
$$ \Delta^{k}([s,t])= \{ (\tau_{1}, \ldots, \tau_{k}) \in [0,T]^{k}: s \leq \tau_{1} \leq \tau_{2} \leq \cdots \leq \tau_{k} \leq t \} $$
for $0\leq s \leq t \leq T$ and $k \in \N$.

By reapplication of the change-of-variable formula (\ref{ito-frac}), 
which holds true for the solution $X^a$ to our SDE,
and setting
$\mathcal{D}^{\alpha} = \mathcal{D}^{\alpha_1} \ldots
\mathcal{D}^{\alpha_k}$ and  $dB^{\alpha}(s, s_1, \ldots,
s_{k-1}) = dB_s^{\alpha_1}dB_{s_1}^{\alpha_2} \ldots$ 
$dB_{s_{k-1}}^{\alpha_{k}}$ for a multi-index $\alpha = (\alpha_1,
\ldots, \alpha_k) \in \{0,1, \ldots, d\}^k$, we get that
\begin{eqnarray} \label{St-proof-expansion-1:Wm}
    f(X_t^a) = f(a) + \sum_{k=1}^{m} \sum_{\alpha \in \{0,1, \ldots, d\}^k}
    \mathcal{D}^{\alpha} f(a) \,\,  \int_{\Delta^{k}([0,t])} \,
    dB^{\alpha}(s, s_1, \ldots, s_{k-1})
    + \mathcal{R}^*_m(0,t)
\end{eqnarray}
with the truncation term
\begin{equation} \label{St-proof-truncation:Wm}
    \mathcal{R}_m^*(0,t) = \sum_{\alpha \in \{0,1, \ldots, d\}^{m+1}}
    \int_{0}^t \int_{0}^{s_{m}} \ldots \int_{0}^{s_1}
    \mathcal{D}^{\alpha} f(X_{s}^a) \, dB_s^{\alpha_1} \, dB_{s_1}^{\alpha_2} \ldots dB_{s_{m}}^{\alpha_{m+1}} .
\end{equation}
Clearly, for $m=0$ there exists only the tree $\textbf{t} = \gamma
\in LTS$ with $l(\textbf{t})-1=0$ and we obtain
\begin{equation}
    \sum_{\substack{\textbf{t} \in LTS \\
        l(\textbf{t})-1 = 0}}
        F(\textbf{t})(a) \,\, \mathcal{I}_{\textbf{t}}(1)_{0,t}
        = F(\gamma)(a) = f(a) .
\end{equation}
Thus, to prove (\ref{Strato-tree-expansion-exact-sol-formula1:Wm})
it is sufficient to show for every $m \in \mathbb{N}$ that
\begin{equation} \label{St-proof-to-show-1}
    \sum_{\alpha \in \{0,1, \ldots, d\}^m} \mathcal{D}^{\alpha} f(a) \,\, \int_{\Delta^{m}([0,t])} \,
    dB^{\alpha}(s, s_1, \ldots, s_{m-1}) =
    \sum_{\substack{\textbf{t} \in LTS \\ l(\textbf{t})-1 = m}}
    F(\textbf{t})(a) \,\, \mathcal{I}_{\textbf{t}}(1)_{0,t} \,\,.
\end{equation}
The proof proceeds by induction. Step $m=1$ is performed for
a better understanding. In this case two different trees
$\textbf{t}_{1}=(\tau_0^2)^1$ and
$\textbf{t}_{2}=(\tau_{j_1}^2)^1$ in $LTS$, all of length 2, with
$\rho(\textbf{t}_1)=1$ and $\rho(\textbf{t}_2)=H$ have to be
considered.
For these two trees we have
\begin{equation}
    \begin{split}
        &\sum_{\alpha \in \{0,1, \ldots, d\}} \mathcal{D}^{\alpha} F(\gamma)(a) \,\,
        \int_{\Delta^{1}([0,t])} \, dB^{\alpha}(s)\\
&\quad        = \sum_{k=1}^n b^k(a) \frac{\partial f}{\partial x^k}(a) \, t
    + \sum_{j_1=1}^d \sum_{k=1}^n \sigma^{k,j_1}(a) \frac{\partial f}{\partial x^k}(a) \,
    \int_0^t \, dB_s^{j_1} \\
        &\quad = F(\textbf{t}_{1})(a) \,\,
        \mathcal{I}_{\textbf{t}_{1}}(1)_{0,t}
        + \sum_{j_1=1}^d F(\textbf{t}_{2})(a) \,\,
        \mathcal{I}_{\textbf{t}_2}(1)_{0,t} \\
        &\quad = \sum_{\substack{\textbf{t} \in LTS \\ l(\textbf{t})-1 = 1}}
    \sum_{j_1, \ldots, j_{s(\textbf{t})}=1}^d
        F(\textbf{t})(a) \,\, \mathcal{I}_{\textbf{t}}(1)_{0,t} \,\,. \notag
    \end{split}
\end{equation}
Under the assumption that equation (\ref{St-proof-to-show-1})
holds for $m \in \mathbb{N}_0$ we proceed to prove the case $m+1$.
Therefore, writing $\alpha=(\alpha_1, \ldots, \alpha_{m+1})$ for
an element of $\{0,1, \ldots, d\}^{m+1}$, we get
\begin{equation} \label{St-proof-eqn9}
\begin{split}
    &\sum_{\alpha \in \{0,1, \ldots, d\}^{m+1}} \mathcal{D}^{\alpha} f(a) \,\,
    \int_{\Delta^{m+1}([0,t])} \,
    dB^{\alpha}(s, s_1, \ldots, s_{m}) \\
    &= \sum_{\substack{\textbf{t} \in LTS \\ l(\textbf{t})-1 = m}}
    \sum_{j_1, \ldots, j_{s(\textbf{t})}=1}^d
    \sum_{\alpha_{m+1} \in \{0,1, \ldots, d\}}
    \mathcal{D}^{\alpha_{m+1}} F(\textbf{t})(a) \,\,
    \int_0^t \mathcal{I}_{\textbf{t}}(1)_{0,s_m} \,
    dB^{\alpha_{m+1}}_{s_m} \,\, .
\end{split}
\end{equation}
Now, we apply Lemma~2.7 in \cite{Roessler} to
$\mathcal{D}^{\alpha_{m+1}} F(\textbf{t})(a)$. Then, it holds for
any $\textbf{u} \in LTS$ with $l(\textbf{u})-1=m$ in the case of
$\alpha_{m+1}=0$ that
\begin{equation} \label{St-proof-eqn7}
    \mathcal{D}^{\alpha_{m+1}} \sum_{j_1, \ldots, j_{s(\textbf{u})}=1}^d
    F(\textbf{u})(a)
    = \sum_{k=1}^n b^k(a) \frac{\partial}{\partial x^k}
    \sum_{j_1, \ldots, j_{s(\textbf{u})}=1}^d
    F(\textbf{u})(a)
    = \sum_{\textbf{t} \in D(\textbf{u})}
    \sum_{j_1, \ldots, j_{s(\textbf{t})}=1}^d
    F(\textbf{t})(a),
\end{equation}
where $D(\textbf{u})$ is the set of all trees $\textbf{t} \in LTS$
with $l(\textbf{t})=m+2$, $\textbf{t}'|_{\{2, \ldots, m+1\} } =
\textbf{u}'$, $\textbf{t}''|_{\{1, \ldots, m+1\} } = \textbf{u}''$
and $\textbf{t}''(m+2) = \tau_0$. Clearly
$s(\textbf{u})=s(\textbf{t})$ holds for all
$\textbf{t} \in D(\textbf{u})$.

\vspace{0.2cm}

Then we proceed by considering the case of $\alpha_{m+1} \in \{1,
\ldots, d\}$. Again, by  applying Lemma~2.7 in \cite{Roessler}, we
get for $\textbf{u} \in LTS$ with $l(\textbf{u})-1=m$ that
\begin{equation} \label{St-proof-eqn5}
    \begin{split}
    \sum_{\alpha_{m+1} \in \{1, \ldots, d\}}
    \mathcal{D}^{\alpha_{m+1}} \sum_{j_1, \ldots, j_{s(\textbf{u})}=1}^d
    F(\textbf{u})(a)
    &= \sum_{j_1, \ldots, j_{s(\textbf{u})}, \alpha_{m+1}=1}^d
    \sum_{k=1}^n \sigma^{k,\alpha_{m+1}}(a) \frac{\partial}{\partial x^k}
    F(\textbf{u})(a) \\
    &= \sum_{\textbf{t} \in S(\textbf{u})} \sum_{j_1, \ldots,
    j_{s(\textbf{t})}=1}^d F(\textbf{t})(a),
    \end{split}
\end{equation}
where $S(\textbf{u})$ denotes the set of trees $\textbf{t} \in
LTS$ with $l(\textbf{t})=m+2$, $\textbf{t}' |_{\{2,\ldots,m+1\}} =
\textbf{u}'$, $\textbf{t}'' |_{\{1, \ldots, m+1\}} = \textbf{u}''$
and $\textbf{t}''(m + 2) = \tau_{j_{s(\textbf{u})+1}}$. Here we have $s(\textbf{t}) =
s(\textbf{u}) + 1$ for all $\textbf{t} \in S(\textbf{u})$. 

\vspace{0.2cm}

Combining now the results for the case of $\alpha_{m+1}=0$
and $\alpha_{m+1} \in \{1, \ldots, d\}$, the equation
\begin{equation} \label{St-proof-eqn10}
\begin{split}
    \sum_{\alpha_{m+1} \in \{0,1, \ldots, d\}} \mathcal{D}^{\alpha^{m+1}}
    \sum_{j_1, \ldots, j_{s(\textbf{u})}=1}^d F(\textbf{u})(a)
    = \sum_{\textbf{t} \in D(\textbf{u}) \cup S(\textbf{u})}
    \sum_{j_1, \ldots, j_{s(\textbf{t})}=1}^d F(\textbf{t})(a)
\end{split}
\end{equation}
holds for every $\textbf{u} \in LTS$ with $l(\textbf{u})-1=m$. Now
it is easily seen that
\begin{equation} \label{St-proof-eqn11}
    \bigcup_{\substack{\textbf{u} \in LTS \\ l(\textbf{u})-1=m}}
    D(\textbf{u}) \cup S(\textbf{u}) = \{ \textbf{t} \in
    LTS : l(\textbf{t})-1 = m+1\}.
\end{equation}
As a last step, we observe that
\begin{equation} \label{St-proof-eqn12}
    \int_0^t \mathcal{I}_{\textbf{u}}(1)_{0,s_m} \,
    dB^{\alpha_{m+1}}_{s_m}
    = \mathcal{I}_{\textbf{t}}(1)_{0,t}
\end{equation}
holds for all $\textbf{t} \in D(\textbf{u})$ in the case of
$\alpha_{m+1}=0$ and for all $\textbf{t} \in S(\textbf{u})$ with
$\textbf{t}''(m+1)=\tau_{\alpha_{m+1}}$ for all $\alpha_{m+1} \in \{1, \ldots, d\}$.
By applying (\ref{St-proof-eqn10})--(\ref{St-proof-eqn12}) to
(\ref{St-proof-eqn9}) we thus arrive at (\ref{St-proof-to-show-1})
with $m$ replaced by $m+1$, which completes the proof of the first
part of Theorem~\ref{St-tree-expansion-exact-sol:Wm}. 

\vspace{0.2cm}

Finally, we have to prove that $\mathcal{R}_m(0,t) =
\mathcal{R}_m^*(0,t)$ given in
(\ref{Strato-tree-expansion-exact-sol-formula1:Wm}) and
(\ref{St-proof-expansion-1:Wm}). From
(\ref{St-proof-eqn7})--(\ref{St-proof-eqn11}) it follows that for
each $\alpha \in \{0,1, \ldots, d\}^{m+1}$ there exists a subset
$LTS(\alpha) \subset LTS$ with a fixed choice of $j_1, \ldots,
j_{s(\textbf{t})} \in \{1, \ldots, d\}$ for $\textbf{t} \in
LTS(\alpha)$ such that $l(\textbf{t})-1=m+1$, $\textbf{t}''(1) =
\gamma$, $\textbf{t}''(i)=\tau_{\alpha_{i-1}}$ for $i=2, \ldots,
m+2$ and
\begin{equation} \label{St-proof-eqn13}
    \mathcal{D}^{\alpha} f(a) = \sum_{\textbf{t} \in LTS(\alpha)}
    F(\textbf{t})(a)
\end{equation}
for all $a \in \mathbb{R}^n$ and $m \in \mathbb{N}_0$. Then, the
sets $LTS(\alpha)$, $\alpha \in \{0,1, \ldots, d\}^{m+1}$, build a
partition of $\{\textbf{t} \in LTS : l(\textbf{t})-1 = m+1\}$.
Thus, we have $\mathcal{I}_{\textbf{t}}(Z_s^a)_{0,t} =
\mathcal{I}_{\textbf{u}}(Z_s^a)_{0,t}$ for all $\textbf{t},
\textbf{u} \in LTS(\alpha)$ and any integrable process $Z$.
Replacing now $Z_s^a$ by $F(\textbf{t})(X_s^a)$ yields that for
all $\alpha \in \{0,1, \ldots, d\}^{m+1}$, the following
relation holds true:
\begin{equation}
    \sum_{\textbf{t} \in LTS(\alpha)}
    \mathcal{I}_{\textbf{t}}(F(\textbf{t})(X_s^a))_{0,t}
    = \int_{0}^t \int_{0}^{s_{m}} \ldots \int_{0}^{s_1}
    \mathcal{D}^{\alpha} f(X_{s}^a) \, dB_s^{\alpha_1}
    \, dB_{s_1}^{\alpha_2} \ldots dB_{s_{m}}^{\alpha_{m+1}}
\end{equation}
which completes the proof.\fin
Invoking Theorem~\ref{St-tree-expansion-exact-sol:Wm} we obtain the
following corollary:
%
%
\begin{corollary} \label{St-tree-expansion-expectation-exact-sol:Wm}
    Let $(X_t^a)_{t \in I}$ be the solution of SDE~(\ref{eqintro})
    with initial value $X_{0}^a = a \in \mathbb{R}^n$.
    Then for $m \in \mathbb{N}_0$ and $f \in  C_b^{m+1}(\mathbb{R}^n, \mathbb{R}), b, \sigma^j \in
    C_b^{m+1}(\mathbb{R}^n, \mathbb{R}^n)$, $1 \leq j \leq d$,
    we get for $t \in [0,T]$ the expansion
    \begin{equation} \label{Strato-tree-expansion-expectation-exact-sol-formula1:Wm}
        {\rm P}_t f(a) = \sum_{\substack{\textbf{t} \in LTS(S) \\
        l(\textbf{t}) \leq m+1}}
        \sum_{j_1, \ldots, j_{s(\textbf{t})}=1}^d
        F(\textbf{t})(a) \,\, {\rm E}(\mathcal{I}_{\textbf{t}}(1)_{0,1})
        \,\, t^{\rho(\textbf{t})} + {\rm E}(\mathcal{R}_m(0,t))
        \,\, .
    \end{equation}
\end{corollary}
\quad \\
\noindent {\it Proof.} Apply Theorem~\ref{St-tree-expansion-exact-sol:Wm}
and take the expectation in formula
(\ref{Strato-tree-expansion-exact-sol-formula1:Wm}). Using the
notation of the proof of
Theorem~\ref{St-tree-expansion-exact-sol:Wm}, we observe that for
each $\textbf{t} \in LTS$ there exists an $\alpha \in \{0,1, \ldots, d\}^{m}$
with $m=l(\textbf{t})-1$ such that $\textbf{t} \in LTS(\alpha)$.
Then, due to (\ref{St-proof-eqn12}) and the scaling property (\ref{scaling})
it follows that
\begin{equation}
    {\rm E} \left( \mathcal{I}_{\textbf{t}}(1)_{0,t} \right) =
    {\rm E} \left( \mathcal{I}_{\textbf{t}}(1)_{0,1} \right)
    t^{H \, |\alpha| + m -|\alpha|}
    = {\rm E} \left(
    \mathcal{I}_{\textbf{t}}(1)_{0,1} \right) \,\,
    t^{\rho(\textbf{t})}
\end{equation}
holds with $|\alpha| = \sum_{i=1}^m 1_{\{\alpha_i \neq 0\}}$ since
we have ${\rho(\textbf{t})} = H \, |\alpha| + m -|\alpha|$.
\fin
%
%
%
%
%

The sequel of the paper is now devoted to derive the announced controls
on the remainder term
${\rm E}({\cal R}_m(0,t))$ appearing in 
(\ref{Strato-tree-expansion-expectation-exact-sol-formula1:Wm}),
according to the value of $H$ and the assumptions on $f,b$ and $\sigma$.
These estimates will imply easily our Theorem \ref{thmptfa}.

\subsection{Study of the remainder term for $1/3<H<1/2$}

We assume in this section that $1/3<H<1/2$ and that assumption
(A) holds true. Then we will show that  for fixed $m\in\N$, we have
\begin{equation}\label{rem13}
{\rm E}({\cal R}_m(0,t))=O(t^{(m+1)H})
\end{equation} as $t\rightarrow 0,$
a fact which trivially yields (\ref{ptfa}) in Theorem \ref{thmptfa}.
Furthermore, notice that the control (\ref{rem13})  is a 
direct consequence of the following:

\begin{lemma}\label{ctrl-lm}
Let $g\in{\rm C}^2(\R^n)$ be bounded together with its derivatives
and $X$ be the unique solution to (\ref{eds-frac}) in $\cq_{\ka,a}(\R^n)$ with
$\ka\in (\frac{1-H}{2},H)$.
For any $\alpha_1,\ldots,\alpha_r\in\{0,\ldots,d\}$ we have
\begin{equation}\label{ctrl-rem}
{\rm E}\left|\int_{\Delta^r([0,t])} g(X)dB^{\alpha_r}\ldots dB^{\alpha_1}\right|=O(t^{r-|\alpha|(1-H)}),\mbox{ as }t\rightarrow 0,
\end{equation}
where $|\alpha| = \sum_{i=1}^r 1_{\{\alpha_i \neq 0\}}$.
\end{lemma}
\noindent{\it Proof}.
The more difficult setting holds when $\alpha_j\neq 0$ for all $1\le j\le r$, that is
when $r=|\alpha|$. For this reason we will only prove the assertion  in this case. Moreover, we split the proof  in several steps.

\vspace{0.3cm}

\noindent
{\it Step 1: Scaling.} For $j\in\{1,\ldots,r\}$ and
$c > 0 $ set $B^{\alpha_j,(c)}_u=c^H\,B^{\alpha_j}_{u/c}$ and let
$X^{(c)}$ denote the solution of (\ref{eds-frac}), where $B$ is replaced by $B^{(c)}$.
For fixed $t$, we have
$$
\begin{array}{lll}
\int_{\Delta^r([0,t])} g(X)dB^{\alpha_r}\ldots dB^{\alpha_1}
&=&\int_0^t dB^{\alpha_1}_{t_1}\int_0^{t_1} dB^{\alpha_2}_{t_2} \ldots\int_0^{t_{r-1}} dB^{\alpha_r}_{t_r}
g(X_{t_r})\\
&=&\int_0^1 dB^{\alpha_1}_{t\cdot t_1}\int_0^{t_1} dB^{\alpha_2}_{t\cdot t_2} \ldots\int_0^{t_{r-1}}
dB^{\alpha_r}_{t\cdot t_r}g(X_{t\cdot t_r})\\
&\stackrel{\mathcal{L}}{=} &\int_0^1 dB^{\alpha_1,(t)}_{t\cdot t_1}\int_0^{t_1} dB^{\alpha_2,(t)}_{t\cdot t_2} \ldots\int_0^{t_{r-1}}
dB^{\alpha_r,(t)}_{t\cdot t_r}g(X^{(t)}_{t\cdot t_r})\\
&=&t^{rH}\int_0^1 dB^{\alpha_1}_{t_1}\int_0^{t_1} dB^{\alpha_2}_{t_2} \ldots\int_0^{t_{r-1}}
dB^{\alpha_r}_{t_r}g(X^{(t)}_{t\cdot t_r}).
\end{array}
$$
Consequently, in order to obtain (\ref{ctrl-rem}), it suffices to prove that
$$\sup_{t\in [0,T]}{\rm E}\left|
\int_0^1 dB^{\alpha_1}_{t_1}\int_0^{t_1} dB^{\alpha_2}_{t_2} \ldots\int_0^{t_{r-1}}
dB^{\alpha_r}_{t_r}\,g(X^{(t)}_{t\cdot t_r})
\right|< \infty.
$$

\vspace{0.3cm}

\noindent
{\it Step 2:} Fix $t$ and set
$$
z_s=\int_0^{s}
dB^{\alpha_r}_{t_r}\,g(X^{(t)}_{t\cdot t_r}),\quad s\in [0,1].
$$
By (\ref{bnd:norm-imdx}) and
(\ref{bnd:phi}) we have
$$
{\cal N}(z,\cq_{\ka,0})\le c_{B}\big(1+{\cal N}(g(X^{(t)}_{t\,\cdot}),\cq_{\ka,g(a)})\big)
\le c_{B}c_{g}\big(1+{\cal N}^2(X^{(t)}_{t\,\cdot},\cq_{\ka,a})\big).
$$
Here, $c_B>1$ is the random constant appearing in (\ref{bnd:norm-imdx}) (see also (\ref{bnd:inc-z})), whose value  will not change from line to line, while
$c_g$ denotes a non-random constant depending only on $g$, whose value  can change from one line to
another.
Set now
$$
q_s=\int_0^{s} dB^{\alpha_{r-1}}_{t_{r-1}}\int_0^{t_{r-1}}
dB^{\alpha_r}_{t_r}g(X^{(t)}_{t\cdot t_r})=\int_0^s dB^{\alpha_{r-1}}_{t_{r-1}}z_{t_{r-1}},\quad s\in [0,1].
$$
Similarly, we have
$$
{\cal N}(q,\cq_{\ka,0})\le c_{B}\big(1+{\cal N}(z,\cq_{\ka,0})\big)
\le {c_{B}}^2c_{g}\big(1+{\cal N}^2(X^{(t)}_{t\,\cdot},\cq_{\ka,a})\big).
$$
By induction, we easily deduce that
$$
{\cal N}\left(\int_0^{\cdot}dB^{\alpha_1}_{t_1}\int_0^{t_1} dB^{\alpha_2}_{t_2} \ldots\int_0^{t_{r-1}}
dB^{\alpha_r}_{t_r}g(X^{(t)}_{t\cdot t_r})
,\cq_{\ka,0}\right)
\le {c_B}^rc_{g}\big(1+
{\cal N}^2(X^{(t)}_{t\,\cdot},\cq_{\ka,a})
\big).
$$
Since we have $|z_1|\le \|z\|_{\ka}\le
{\cal N}(z,\cq_{\ka,0})$ for a path $z$ starting from $0$, we deduce from the Cauchy-Schwarz inequality
that (\ref{ctrl-rem}) is in fact a consequence of
showing
\begin{equation}\label{cst}
{\rm E}({c_B}^{2r})<\infty
\end{equation}
and
\begin{equation}\label{n4}
\sup_{t\in[0,T]}{\rm E}\left|
{\cal N}^4(X^{(t)}_{t\,\cdot},\cq_{\ka,a})\right|
=\sup_{t\in[0,T]}{\rm E}\left|
{\cal N}^4(X_{t\,\cdot},\cq_{\ka,a})\right|
< \infty.
\end{equation}

\vspace{0.3cm}

\noindent
{\it Step 3:} Using that
$B$ has moments of all order and (\ref{lem:garsia}), we easily obtain by (\ref{bnd:inc-z}) that
(\ref{cst}) is verified.
So, let us concentrate on (\ref{n4}), which is more difficult.
We will in fact only prove that ${\rm E}
\big|{\cal N}^4(X,\cq_{\ka,a})\big|<+\infty$, since we can obtain the
control (\ref{n4}) in a
similar way. Recall from the proof of Theorem \ref{thm:ex-uniq} that $X$ defined on $[0,\tau]$
belongs by definition to the ball $B_M$ given by (\ref{def:ball-m}), where $M$ and $\tau$ verify
$$
M\ge c_{\sigma,B}\big(1+\tau^{\gamma-\kappa}M^2\big).
$$
For fixed $\tau$, the inequality $u\ge c_{\sigma,B}\big(1+\tau^{\gamma-\kappa}u^2\big)$
admits solutions $u$ iff ${c_{\sigma,B}}^{-2}-4\tau^{\gamma-\ka}>0$, i.e., iff $\tau^{\gamma-\ka}<(4\,{c_{\sigma,B}}^2)^{-1}$. In this case,
the solutions are $u\in [M_-,M_+]$, whereby
$$
M_\pm=\frac{
\frac{1}{c_{\sigma,B}}
\pm\sqrt{
\frac{1}{{c_{\sigma,B}}^2}
-4\tau^{\gamma-\ka}
}
}
{2\tau^{\gamma-\ka}}.
$$
By choosing for instance 
\begin{equation}\label{tau-1}
\tau^{\gamma-\ka}=(8\,{c_{\sigma,B}}^2)^{-1}
\end{equation} 
we obtain that
\begin{equation}\label{tau-2}
{\cal N}(X_{|_{[0,\tau]}},\cq_{\ka,a})
\le (4-2\sqrt{2})c_{\sigma,B}.
\end{equation} 
Furthermore,
as explained at the end of the proof of Theorem \ref{thm:ex-uniq} and due to the crucial fact
that $\sigma$ and  its derivatives are bounded, we can
in fact choose the same $M$ for the bound of $\delta X$ on $[\tau,2\tau]$, $[2\tau,3\tau]$, etc.
Using the triangle inequality we deduce:
$$
\begin{array}{lll}
&{\cal N}(X,\cq_{\ka,a})\\
& \quad \le{\cal N}(X_{|_{[0,\tau]}},\cq_{\ka,a})+{\cal N}(X_{|_{[\tau,2\tau]}},\cq_{\ka,X_\tau})+\ldots
+{\cal N}(X_{|_{[\lfloor T\tau^{-1}\rfloor\tau,T]}},\cq_{\ka,X_{\lfloor T\tau^{-1}\rfloor\tau}})\\
& \quad \le(\lfloor T\tau^{-1}\rfloor+1)M.
\end{array}
$$
In other words, we deduce ${\cal N}(X,\cq_{\ka,a})\le  {\rm cst}\,{c_{\sigma,B}}^{1+\frac{2}{\gamma-\kappa}}$,
see (\ref{tau-1}) and (\ref{tau-2}).
Thus it follows easily  that ${\rm E}
\big|{\cal N}^4(X,\cq_{\ka,a})\big|<+\infty$ and the proof of Lemma \ref{ctrl-lm} is finished.
\fin

\bigskip

\subsection{Some properties of iterated integrals in the case $H>1/2$}

Let us say first a few words about the strategy we have adopted
in order to get  equation (\ref{ptfa}): the key point will be
again to get an accurate bound for $E[{\mathcal R}_m(0,t)]$, and
thus we use estimates based on Malliavin
calculus tools and explicit computations of moments for multiple iterated
integrals with respect to the fractional Brownian motion. For a proposed
pathwise
control on the remainder ${\mathcal R}_m(0,t)$ , see.,   e.g. in \cite[Remark 7.4]{Gu2}.

Before we turn to the control of the remainder in the case $H>1/2$, we will establish first some properties of  iterated integrals with respect to fractional Brownian motion. To do this, we require some additional notations.

For a multi-index $\alpha \in \{0,1, \ldots , d\}^{k}$ with $k \in \N$  denote by $l(\alpha)$ the
length of $\alpha$, i.e., $l(\alpha)=k$.
Moreover set $\mathcal{A}_k=\{0,1,\ldots,d\}^k$ for $k\in\N$, i.e.,
$\mathcal{A}_{k}$ is the set of all multi-indices of length $k$.
Furthermore, define for $\alpha \in \mathcal{A}_{k}$ the sets
$$ \mathfrak{J}_{\alpha}= \{ j=1, \ldots, k: \alpha_{j} \neq 0 \},
\quad\mbox{ and }\quad
\mathfrak{J}_{\alpha,i }= \{ j=1, \ldots, k: \alpha_{j} =i  \},
$$
for $i=1, \ldots, d$ and
$|\alpha|  = |\mathfrak{J}_{\alpha}| .$
Finally for a multi-index $\alpha\in \mathcal{A}_{k}$
and $j=1, \ldots, k$ we  denote
$$ \alpha_{-j}=(\alpha_{1},\alpha_{2}, \ldots, \alpha_{j-1},\alpha_{j+1},\ldots, \alpha_{k}).$$
Recall that for $m \in \N$ and $0 \leq t_{1} \leq t_{2} \leq T$ we set
$$ \Delta^{m}([t_{1},t_{2}])= \{ (\tau_{1}, \ldots, \tau_{m}) \in [0,T]^{m}: t_{1} \leq \tau_{1} \leq \tau_{2} \leq \cdots \leq \tau_{m} \leq t_{2} \}.$$
Moreover, we will use the notation
$$\int_{\Delta^{k}([t_{1},t_{2}])} dB^{\alpha} = \int_{t_{1}}^{t_{2}} \int_{t_{1}}^{s_{k-1}} \cdots \int_{t_{1}}^{s_1}
     \, dB_s^{\alpha_1} \, dB_{s_1}^{\alpha_2} \ldots dB_{s_{k}}^{\alpha_{k}} $$
for $\alpha \in \mathcal{A}_{k}$.

With these notations in hand, the following proposition is shown
easily, and its proof will be omitted here. Indeed,
part (a) follows immediately by the symmetry
of fractional Brownian motion and part (b) can be shown analogously to 
Theorem 11 in \cite{BC}.

\bigskip

\begin{proposition}\label{mom_12}
Let $k \in \N$ and  $\alpha \in \mathcal{A}_{k}$.  \\
(a) If $|\alpha|$ is odd, then we have
$$  {\rm E}\int_{\Delta^{k}([0,1])} dB^{\alpha} =0.$$
(b) If $|\alpha|$ is even, then it holds
$${\rm E}\int_{\Delta^{k}([0,1])} dB^{\alpha} = \frac{(\gamma_{H}/2)^{|\alpha|/2}}{(|\alpha|/2)! } \sum_{\mathfrak{s} \in \mathfrak{S}_{\mathfrak{J}_{\alpha}}} \mathcal{V}(\mathfrak{s}, \alpha)$$
with
$$ \mathcal{V}(\mathfrak{s}, \alpha)= \int_{0\le t_1<\ldots<t_k\le 1}
 \prod_{l=1}^{|\alpha|/2} \delta_{\alpha_{\mathfrak{s}(2l-1)},\alpha_{\mathfrak{s}(2l)}} |t_{\mathfrak{s}(2l)}-t_{\mathfrak{s}(2l-1)}|^{2H-2} dt_{1} \ldots  dt_{k} ,   $$
where $\mathfrak{S}_{\mathfrak{J}_{\alpha}}$ is the group of all permutations of the set $\mathfrak{J}_{\alpha}$, $\gamma_{H}=H(2H-1)$ and  $\delta_{i,j}$ is Kronecker's symbol. \\
(c) It holds
\begin{eqnarray*}
&&   {\rm E} \left| \int_{\Delta^{k}([0,1])} dB^{\alpha}\right|^{2} = \frac{(\gamma_{H}/2)^{|\alpha| }}{|\alpha|! }   \sum_{\mathfrak{s} \in \mathfrak{S}^{\mathfrak{J}_{\alpha}^{2}}}   \mathcal{W}(\mathfrak{s}, \alpha)
\end{eqnarray*}
with
\begin{align*}
 & \mathcal{W}(\mathfrak{s}, \alpha)= \\ & \,\, 
\int_{0\le t_1<\ldots<t_k\le 1} \int_{0\le t_{k+1}<\ldots<t_{2k}\le 1}   \prod_{l=1}^{|\alpha|} \delta_{\alpha_{\mathfrak{s}(2l-1)},\alpha_{\mathfrak{s}(2l)}}  |t_{\mathfrak{s}(2l)}-t_{\mathfrak{s}(2l-1)}|^{2H-2} \, d t_{k+1} \cdots\, d t_{2k} d t_{1} \cdots \, d t_{k} ,\end{align*}
where $\mathfrak{S}^{\mathfrak{J}_{\alpha}^{2}}$ denotes the group of all permutations of the set $$\mathfrak{J}_{\alpha}^{2}= \{j=1, \ldots, 2k: j \in \mathfrak{J}_{\alpha} \textrm{ or } j-k  \in  \mathfrak{J}_{\alpha}  \}.$$

\end{proposition}
Notice that part (a) and (b) of the above proposition yield a representation  for  the coefficients ${\rm E}(\mathcal{I}_{\textbf{t}}(1)_{0,1})$, since $$ {\rm E}(\mathcal{I}_{\textbf{t}}(1)_{0,1})=  {\rm E} \int_{\Delta^{l(\bf{t})}([0,1])} d B^{\alpha} 
$$ with ${\bf t''}(i)=\alpha_{i} \in \{0,1, \ldots, d\}$, $i=1, \ldots, l({\bf t})$. 

\bigskip

For further computations, we also need  the following positivity result
for iterated integrals of the fractional Brownian motion.

\bigskip

\begin{proposition}\label{mix_mom} Let  $m_{i} \in \N$ for $i=1, \ldots, n$ with $n \in \N$. Moreover let $\alpha^{m_{i}} \in \mathcal{A}_{m_{i}}$ and $0\leq s_{i} \leq t_{i} \leq T$ for  $i=1, \ldots, n$. It holds
$$ {\rm E} \left[  \prod_{i=1}^{n} \int_{\Delta^{m_{i}}([s_{i},t_{i}])} dB^{\alpha^{m_{i}}} \right]  \geq 0.$$
\end{proposition}

\noindent{\it Proof.}
Let $t_{k}^{l}=k2^{-l}$, $k=0, 1, \ldots, 2^{l}$.
Denote by $B^{l,(\alpha_{j}^{m_{i}})}$ the piecewise linear interpolation of $B^{(\alpha_{j}^{m_{i}})}$ with step size $2^{-l}$, i.e.,
$$B_{t}^{l,(\alpha_{j}^{m_{i}})}=B_{t_{k}^{l}}^{(\alpha_{j}^{m_{i}})} + 2^{l}(t-t_{k}^{l})(B_{t_{k+1}^{l}}^{(\alpha_{j}^{m_{i}})}-B_{t_{k}^{l}}^{(\alpha_{j}^{m_{i}})}), \qquad t \in [t_{k}^{l},t_{k+1}^{l}).$$
We have
$$ \prod_{i=1}^{n} \int_{\Delta^{m_{i}}([s_{i},t_{i}])}dB^{\alpha^{m_{i}}} = \lim_{l \rightarrow \infty} \prod_{i=1}^{n}  \int_{\Delta^{m_{i}}([s_{i},t_{i}])}dB^{l,\alpha^{m_{i}}}  $$
almost surely, due to Proposition 3.12.
Since $ \prod_{i=1}^{n} \int_{\Delta^{m_{i}}([s_{i},t_{i}])}dB^{l,\alpha^{m_{i}}}$ belongs to a finite Wiener chaos, we also have 
$$ {\rm E} \prod_{i=1}^{n} \int_{\Delta^{m_{i}}([s_{i},t_{i}])}dB^{\alpha^{m_{i}}} =  \lim_{l \rightarrow \infty}  {\rm E} \prod_{i=1}^{n}\int_{\Delta^{m_{i}}([s_{i},t_{i}])}dB^{l,\alpha^{m_{i}}}   $$ 
according to \cite{B}.  Note that
$$  \int_{\Delta^{m_{i}}([s_{i},t_{i}])}dB^{l,\alpha^{m_{i}}} =  \int_{s_{i}}^{t_{i}} \int_{s_{i}}^{t_{m_{i}}} \cdots \int_{s_{i}}^{t_{2}} \prod_{j=1}^{m_{i}} Z_{t_j}^{l,(\alpha_{j}^{m_{i}})}  \, dt_{1} \cdots \, d t_{m_{i}-1}\, d t_{m_{i}}, $$
where
$$ Z_{t}^{l,(\alpha_{j}^{m_{i}})}=  2^{l}(B^{(\alpha_{j}^{m_{i}})}_{t_{k+1}^{l}} -B^{(\alpha_{j}^{m_{i}})}_{t_{k}^{l}} ), \qquad t \in [t_{k}^{l},t_{k+1}^{l}) $$
for $\alpha_{j}^{m_{i}}\neq 0$ and  $Z_{t}^{l,(0)}=1$ for $t \in [0,T]$. We thus  have
\begin{eqnarray*}
&& \prod_{i=1}^{n} \int_{\Delta^{m_{i}}([s_{i},t_{i}])}dB^{l,\alpha^{m_{i}}} \\ && \qquad =
\int_{\Delta^{m_{n}}([s_{n},t_{n}])}  \cdots \int_{\Delta^{m_{1}}([s_{1},t_{1}])}  \prod_{i=1}^{n} \prod_{j=1}^{m_{i}} Z_{t_j^{i}}^{l,(\alpha_{j}^{m_{i}})}  dt_{1}^{m_{1}} \ldots d t_{m_{1}}^{m_{1}} \cdots d t_{1}^{m_{n}} \ldots d t_{m_{n}}^{m_{n}}.
\end{eqnarray*}
But the term   $\prod_{i=1}^{n} \prod_{j=1}^{m_{i}} Z_{t_j^{i}}^{l,(\alpha_{j}^{m_{i}})} $ is a 
product, which consists only of increments of the independent fractional Brownian motions 
$B^{1}, \ldots, B^{d}$ 
with Hurst parameter $H>1/2$ and of the constant factors 1. Since it is well-known that
the increments of a fractional Brownian motion of Hurst index $H>1/2$ are positively
correlated, and also that we have, for a centered Gaussian vector $(G_1,\ldots,G_{2k})$:
$$
{\rm E}(G_1\ldots G_{2k})=\frac{1}{k!2^k}\sum_{\mathfrak{s}\in\mathfrak{S}_{2k}}\prod_{\ell=1}^k
{\rm E}(G_{\mathfrak{s}(2\ell)}G_{\mathfrak{s}(2\ell-1)}),
$$
we clearly deduce that
$$  {\rm E} \, \prod_{i=1}^{n} \prod_{j=1}^{m_{i}} Z_{t_j^{i}}^{l,(\alpha_{j}^{m_{i}})} \geq 0$$
for all $t_{1}^{m_{1}}, \ldots ,t_{m_{n}}^{m_{n}} \in [0,T]$.
Hence we obtain 
$$  {\rm E} \prod_{i=1}^{n}\int_{\Delta^{m_{i}}([s_{i},t_{i}])}dB^{l,\alpha^{m_{i}}}   \geq  0 $$
for every $l \in \N$, and the assertion follows.
\begin{flushright} $\square$ \end{flushright}

Our estimate of the remainder will
also require the Malliavin derivative of an iterated integral.
Recall then that for a random variable $F\in \sk^{1,2}$ we denote by
$D^{i}F$ the $i$-th component of the Malliavin derivative, i.e., 
$D F = (D^{1}F, \ldots, D^{d}F).$
Recall moreover that for $\alpha \in \mathcal{A}_{k}$, $k \in \N$, we have defined
$$ \mathfrak{J}_{\alpha}= \{ j=1, \ldots, k: \alpha_{j} \neq 0 \},
\qquad
\mathfrak{J}_{\alpha,i }= \{ j=1, \ldots, k: \alpha_{j} =i  \}, $$
for $i=1, \ldots, d$ 
and 
$$ \alpha_{-j}=(\alpha_{1},\alpha_{2}, \ldots, \alpha_{j-1},\alpha_{j+1},\ldots, \alpha_{k}).$$
Then the stochastic derivative of a multiple integral can be computed
as follows:
\begin{proposition}\label{mal_deriv}
Let $m \in \N$ and $\alpha \in \mathcal{A}_{m}$. We have
\begin{eqnarray}{\label{deriv}}
 &&  D_{u}^{i} \int_{\Delta^{m}([s,t])} d B^{\alpha}(t_{1}, \ldots, t_{m})\, \\&& \,\,\qquad \qquad \qquad  = \sum_{j \in \mathfrak{J}_{\alpha,i}}\int_{   s \leq  t_{1} \leq  \ldots \leq t_{j-1}\leq u \leq t_{j} \leq \ldots \leq t_{m-1}\leq t } \, d B^{\alpha_{-j}}(t_{1}, \ldots, t_{m-1})  \nonumber
 \end{eqnarray}  for $i=1, \ldots, d$.
\end{proposition}

\bigskip
\noindent {\it Proof.}
We proceed by induction over $l(\alpha)$. 

\vspace{0.5cm}

\noindent
(a) Assume that $l(\alpha)=1$. For  $\alpha=(0)$ the assertion clearly holds.
Moreover for $\alpha=(j)$, $j=1, \ldots, d$, we have
\begin{eqnarray*}
 \ D_{u}^{i} \int_{s}^{t} \, d B_{\tau}^{(j)}   =   D_{u}^{i}  (B_{t}^{(j)}-B_{s}^{(j)})   =  \delta_{i,j}1_{[s,t]}(u),
\quad\mbox{ for }\quad
i=1, \ldots, d,
\end{eqnarray*} 
which corresponds to expression (\ref{deriv}).

\vspace{0.5cm}

\noindent
(b) Now assume that (\ref{deriv}) holds for all multi-indices of length $m$ and all $i=1, \ldots, d$.
For $\alpha \in \mathcal{A}_{m+1}$ we have
$$   D_{u}^{i} \int_{\Delta^{m+1}([s,t])} d B^{\alpha} = D_{u}^{i} \int_{s}^{t} Y_{\tau} \, dB^{(\alpha_{m+1})}_{\tau}$$ with  $$Y_{\tau} =\int_{\Delta^{m}([s,\tau])} d B^{\widetilde{\alpha}},
 $$
where
$\widetilde{\alpha}=\alpha_{-(m+1)}$.
If $\alpha_{m+1}=0$, then
$$    D_{u}^{i} \int_{s}^{t} Y_{\tau} \, d\tau =   \int_{s}^{t} D_{u}^{i} Y_{\tau}\, d\tau .$$ Hence we obtain by the induction assumption
\begin{eqnarray*}
&&  D_{u}^{i} \int_{\Delta^{m+1}([s,t])} d B^{\alpha}(t_{1}, \ldots, t_{m+1})
\\ && \qquad \qquad =\sum_{j \in \mathfrak{J}_{\widetilde{\alpha},i}} \int_{s}^{t}
 \int_{  s \leq  t_{1} \leq  \ldots \leq t_{j-1}\leq u \leq t_{j} \leq \ldots \leq t_{m-1}\leq \tau}
 d B^{\widetilde{\alpha}_{-j}}(t_{1}, \ldots, t_{m-1})  \, d \tau
\\ &&  \qquad \qquad =\sum_{j \in \mathfrak{J}_{\alpha,i}}
 \int_{  s \leq  t_{1} \leq  \ldots \leq t_{j-1}\leq u \leq t_{j} \leq \ldots \leq \tau \leq t }
 d B^{\alpha_{-j}}(t_{1}, \ldots, t_{m-1}, \tau )
,
\end{eqnarray*}
which shows the assertion in this case. Hence it remains to consider the case $\alpha_{m+1} \neq 0$.
In this case, some standard arguments based on the linear interpolation of
$Y$  and Lemma 1.2.3 in \cite{nual} yield
$$   D_{u}^{i}\int_{s}^{t} Y_{\tau} \, d B_{\tau}^{(\alpha_{m+1})}  = \int_{s}^{t} D_{u}^{i} Y_{\tau} \, d B_{\tau}^{(\alpha_{m+1})} + Y_{u}\delta_{i,\alpha_{m+1}}  1_{[s,t]}(u),
$$ for  $i=1, \ldots, d$.
But now we obtain   by the induction assumption that
\begin{eqnarray*} \hspace*{-0.5cm}
&&  D_{u}^{i} \int_{s}^{t} Y_{\tau} \, d B_{\tau}^{(\alpha_{m+1})}
 \\  && \qquad =   \int_{s}^{t} \sum_{j \in \mathfrak{J}_{\widetilde{\alpha},i}}\int_{ s \leq  t_{1} \leq  \ldots \leq t_{j-1}\leq u \leq t_{j} \leq \ldots \leq t_{m-1}\leq \tau }  \, d B^{\widetilde{\alpha}_{-j}}(t_{1}, \ldots, t_{m-1}) \, d B_{\tau}^{(\alpha_{m+1})}\\ && \qquad  \qquad +    \delta_{i,\alpha_{m+1}}  1_{[s,t]}(u)\int_{\Delta^{m}([s,u])} \, d B^{\alpha_{-(m+1)}}  \\
 &&  \qquad  =   \sum_{j \in \mathfrak{J}_{\alpha, i}, j \neq m+1 }\int_{   s \leq  t_{1} \leq  \ldots \leq t_{j-1}\leq u \leq t_{j} \leq \ldots \leq t_{m-1}\leq \tau \leq t }  \, d B^{\alpha_{-j}}(t_{1}, \ldots, t_{m-1}, \tau )  \\  && \qquad   \qquad   +
\delta_{i,\alpha_{m+1}} 1_{[s,t]}(u) \int_{\Delta^{m}([s,u])} \, d B^{\alpha_{-(m+1)}}
\\ && \qquad  
 =  \sum_{j \in \mathfrak{J}_{\alpha,i} }\int_{ s \leq  t_{1} \leq  \ldots \leq t_{j-1}\leq u \leq t_{j} \leq \ldots \leq t_{m-1}\leq t_{m} \leq t } \, d B^{\alpha_{-j}}(t_{1}, \ldots, t_{m-1}, t_{m}),
\end{eqnarray*}
which is our announced relation.
\begin{flushright} $\square$ \end{flushright}

\bigskip

Now, we will establish an estimate for the second moment of an iterated integral, which will be the key for the control of the remainder 
$\mathcal R_m(0,t)$ in the expansion
of $P_tf(a)$. Indeed, the term $(m!)^{-1/2}$ appearing in (\ref{step1})
will be crucial in order to get some series convergence which will 
entail a nice bound on $\mathcal R_m(0,t)$.

\bigskip

\begin{proposition} \label{est_int} Let $m \in \N$ and $\alpha \in \mathcal{A}_{m}$. There exists a constant $K_{1}>0$, depending only on $H$ and $T$, such that
\begin{eqnarray}\label{step1}
&&   \left( {\rm E} \left| \int_{\Delta^{m}([s,t])} dB^{\alpha}\right|^{2} \right)^{1/2}  \leq   \frac{K^{m}_{1}}{ \sqrt{m!}} |t-s|^{|\alpha|H+m-|\alpha|}
\end{eqnarray}
for $0 \leq s \leq t \leq T$.

\end{proposition}

\noindent{\it Proof.} The proof is separated in three steps.

\vspace{0.3cm}

\noindent
(i) By stationarity (\ref{stationarity}) of the fractional Brownian motion, it follows
$$    \,  \int_{\Delta^{m}([s,t])} dB^{\alpha} \stackrel{\mathcal{L}}{=}     \int_{\Delta^{m}([0,t-s])} dB^{\alpha}. $$
Hence we obtain by the scaling property (\ref{scaling}) of fractional Brownian motion that
$$   \,  \int_{\Delta^{m}([s,t])} dB^{\alpha} \stackrel{\mathcal{L}}{=}    (t-s)^{H|\alpha| + m-|\alpha| } \int_{\Delta^{m}([0,1])} dB^{\alpha}   .$$
Now from Proposition \ref{mom_12} (c) it is obvious that we have 
\begin{eqnarray} \label{help_est_0}
 {\rm E} \left| \int_{\Delta^{m}([0,1])} dB^{\alpha}\right|^{2} \leq {\rm E} \left| \int_{\Delta^{m}([0,1])} dB^{\widetilde{\alpha}}\right|^{2}, \end{eqnarray}
where ${\widetilde{\alpha}}$ is given by 
$ \widetilde{\alpha}=  (\widetilde{\alpha}_{1}, \ldots ,  \widetilde{\alpha}_{m})$ with $ \widetilde{\alpha}_{j}=0$ if $j \in \mathfrak{J}_{\alpha,0}$ and   $ \widetilde{\alpha}_{j}=1$ if $j \in \mathfrak{J}_{\alpha}$, i.e., 
all integrals with respect to $B^{(i)}$, $i=2, \ldots, n$ are replaced by integrals with respect to $B^{(1)}$.

\vspace{0.3cm}

\noindent
(ii) In the next step, we will replace also the integrals with respect to $t$ by integrals with respect to $B^{(1)}$. More precisely, we will show that
\begin{eqnarray}\label{help_est}
{\rm E} \left| \int_{\Delta^{m}([0,1])} dB^{\widetilde{\alpha}}\right|^{2} \leq \gamma_{H}^{|\alpha|-m} {\rm E} \left| \int_{\Delta^{m}([0,1])} dB^{(1, \ldots, 1)} \right|^2,
\end{eqnarray}
with $\gamma_H=H(2H-1)$.
To prove (\ref{help_est}) assume first that there is only one integral with respect to $t$, i.e. $|\mathfrak{J}_{\alpha,0}|=1$.
Thus we have
$$  \int_{\Delta^{m}([0,1])} dB^{\widetilde{\alpha}} =   \int_{\Delta^{k_{1}}([0,1])}  \int_{0}^{s} \int_{\Delta^{k_{2}}([0,s])} dB^{\widetilde{\alpha}_{2}} \,ds  \, dB^{\widetilde{\alpha}_{1}}$$
with
$k_{1}+k_{2}+1=m $ and $\widetilde{\alpha}=(\widetilde{\alpha}_{2},0,\widetilde{\alpha}_{1}). $
By rearranging the order of integration, which is possible, since all integrals are pathwise defined, we get
\begin{eqnarray*}
 \int_{\Delta^{k_{1}}([0,1])}  \int_{0}^{s} 
\int_{\Delta^{k_{2}}([0,s])} dB^{\widetilde{\alpha}_{2}} \,ds  \, dB^{\widetilde{\alpha}_{1}}
=
\int_{0}^{1}   Y_s \,ds, 
\end{eqnarray*}
where we have set
$$
Y_{s} =\int_{\Delta^{k_{1}}([s,1])}   \int_{\Delta^{k_{2}}([0,s])} dB^{\widetilde{\alpha}_{2}} \, dB^{\widetilde{\alpha}_{1}}.
$$
With this notation in hand, observe that we also have
$$
\int_{\Delta^{m}([0,1])} dB^{(1, \ldots, 1)}
=
\int_{0}^{1}   Y_s dB^{(1)}_s.
$$
Hence, when $|\mathfrak{J}_{\alpha,0}|=1$, one can recast (\ref{help_est}) into
\begin{equation}\label{help-est2}
{\rm E} \left| \iou Y_s \, ds \right|^{2} \leq \gamma_{H}
{\rm E} \left| \int_{0}^{1}   Y_s dB^{(1)}_s \right|^2.
\end{equation}
We will now proceed to the estimation of the two terms in (\ref{help-est2}):
first of all, we easily get
\begin{eqnarray*}
&& {\rm E}
 \left| \int_{0}^{1}Y_{s} \,ds \right|^{2}   =   \int_{0}^{1} \int_{0}^{1} {\rm E} Y_{s_{1}} Y_{s_{2}}   \,ds_{1} \, ds_{2}.
\end{eqnarray*} 
Let us compute now ${\rm E} | \int_{0}^{1}Y_{s}dB^{(1)}(s) |^{2}$:
by the relation  between the Young and the divergence integral for fractional Brownian motion, see, e.g. \cite{AMN} or Proposition 5.2.3 in \cite{nual},  we have
$$   \int_{0}^{1}Y_{s}dB^{(1)}(s)=  \delta^{(1)}(Y 1_{[0,1]})+     \gamma_{H} \int_{0}^{1} \int_{0}^{1} D_{s_{1}}^{1}Y_{s_{2}} |s_{1}-s_{2}|^{2H-2} \, ds_{1} \, ds_{2},$$
where we use the  notation
$$ \delta^{(1)}(Y 1_{[0,1]})  =  \delta \left((Y 1_{[0,1]}, \ldots, 0) \right).     $$
Thus we obtain 
\begin{eqnarray*}
 && {\rm E} \left| \int_{0}^{1}Y_{s}dB^{(1)}(s) \right|^{2} =   {\rm E} \left|\delta^{(1)}(Y 1_{[0,1]})  \right|^{2}  +  \gamma_{H}^{2} {\rm E} \left| \int_{0}^{1} \int_{0}^{1} D_{s_{1}}^{1}Y_{s_{2}} |s_{1}-s_{2}|^{2H-2} \, ds_{1} \, ds_{2} \right|^{2} \\ && \qquad \qquad \qquad \qquad  \qquad  \quad + 2  \gamma_{H}  {\rm E}  \delta^{(1)}(Y 1_{[0,1]})    \int_{0}^{1} \int_{0}^{1} D_{s_{1}}^{1}Y_{s_{2}} |s_{1}-s_{2}|^{2H-2} \, ds_{1} \, ds_{2} \\ && \qquad \,\,  \quad \quad \qquad \qquad \geq 
  {\rm E} \left|\delta^{(1)}(Y 1_{[0,1]})  \right|^{2}  \\ && \qquad \qquad \qquad \qquad  \qquad  \quad + 2  \gamma_{H}  {\rm E}  \delta^{(1)}(Y 1_{[0,1]})    \int_{0}^{1} \int_{0}^{1} D_{s_{1}}^{1}Y_{s_{2}} |s_{1}-s_{2}|^{2H-2} \, ds_{1} \, ds_{2}
.
\end{eqnarray*}
Since clearly $ \int_{0}^{1} \int_{0}^{1} D_{s_{1}}^{1}Y_{s_{2}} |s_{1}-s_{2}|^{2H-2} \, ds_{1} \, ds_{2}\in \sk^{1,2}$, we have,
owing to (\ref{eq:duality}), that
\begin{eqnarray*}
&& {\rm E} \lc \delta^{(1)}(Y 1_{[0,1]})    \int_{0}^{1} \int_{0}^{1} D_{s_{1}}^{1}Y_{s_{2}} |s_{1}-s_{2}|^{2H-2} \, ds_{1} \, ds_{2} \rc\\ 
&&  \qquad  = \gamma_{H} \int_{0}^{1} \int_{0}^{1}
\int_{0}^{1} \int_{0}^{1} 
{\rm E}\lc Y_{s_{1}} D_{s_2}^{1}  D_{s_{3}}^{1} Y_{s_{4}}\rc 
|s_{3}-s_{4}|^{2H-2} |s_{1}-s_{2}|^{2H-2}  
\, ds_{3} \, ds_{4}\, ds_{1} \, ds_{2}.
\end{eqnarray*} 
By  the definition of $Y_{s}, s \in [0,1]$, and applying Proposition \ref{mal_deriv}, we can decompose the product $Y_{s_{1}} D_{s_{2}}^{1} D_{s_3}^{1}Y_{s_{4}}$ into a  sum of  products of iterated integrals, and hence
$${ \rm E} \lc Y_{s_{1}} D_{s_{2}}^{1} D_{s_3}^{1}Y_{s_{4}}\rc  \geq 0, \qquad s_{1},s_{2},s_{3},s_{4} \in [0,1],$$ by Proposition \ref{mix_mom}. 
Consequently we obtain
 \begin{eqnarray*}
 && {\rm E} \left| \int_{0}^{1}Y_{s}dB^{(1)}(s) \right|^{2} \geq   {\rm E} \left| \delta^{(1)}(Y 1_{[0,1]}) \right|^{2}.
\end{eqnarray*}
Furthermore, invoking \cite{AMN}, we get 
\begin{multline*}  
{\rm E} \left| \delta^{(1)}(Y 1_{[0,1]}) \right|^{2} =  \gamma_{H}  \int_{0}^{1} \int_{0}^{1} {\rm E}  Y_{s_{1}} Y_{s_{2}} |s_{1}-s_{2}|^{2H-2}  \,ds_{1} \, ds_{2} \\ 
  +   \gamma_{H}^{2}\int_{0}^{1} \int_{0}^{1} \int_{0}^{1} \int_{0}^{1} {\rm E}   D_{\tau_{1}}^{1} Y_{s_{1}}  D_{s_{2}}^{1} Y_{\tau_{2}} |s_{1}-s_{2}|^{2H-2}  |\tau_{1}-\tau_{2}|^{2H-2}  \,ds_{1} \, ds_{2}  \,d\tau_{1} \, d\tau_{2}. 
\end{multline*}
Besides, according to Proposition \ref{mix_mom}, and thanks to the fact
that both $Y_{s_{1}} Y_{s_{2}}$ and 
$D_{\tau_{1}}^{1} Y_{s_{1}}  D_{s_{2}}^{1} Y_{\tau_{2}}$
are products of iterated integrals, we obtain that
$$ 
{\rm E}  Y_{s_{1}} Y_{s_{2}} \geq 0, 
\quad\mbox{ and }\quad
{\rm E} D_{\tau_{1}}^{1} Y_{s_{1}}  D_{s_{2}}^{1} Y_{\tau_{2}} \geq 0,
\qquad s_{1},s_{2} \in [0,1]. 
$$
Since 
$$ |s_{1}-s_{2}|^{2H-2} \geq 1 , \qquad s_{1},s_{2} \in [0,1],$$
we end up with
$$  {\rm E}
 \left| \int_{0}^{1}Y_{s} \,dB^{(1)}(s) \right|^{2} \geq \gamma_{H} {\rm E}
 \left| \int_{0}^{1}Y_{s} \,d s \right|^{2}, $$
which is the announced relation (\ref{help-est2}).
We have thus proved that
\begin{multline*} 
{\rm E} \left| \int_{\Delta^{k_{1}}([0,1])}  
\int_{0}^{s} \int_{\Delta^{k_{2}}([0,s])} dB^{\widetilde{\alpha}_{2}} \,ds  
\, dB^{\widetilde{\alpha}_{1}} \right|^{2}\\ 
\leq \gamma_{H}^{-1}  {\rm E} \left| \int_{\Delta^{k_{1}}([0,1])}  \int_{0}^{s} \int_{\Delta^{k_{2}}([0,s])} dB^{\widetilde{\alpha}_{2}} \,dB^{(1)}  \, dB^{\widetilde{\alpha}_{1}} \right|^{2}.
\end{multline*}
Applying this procedure $m-|\alpha|$ times to replace all integrals with respect to $t$,  equation (\ref{help_est}) is now easily checked. 

\vspace{0.3cm}

\noindent
(iii) Let us conclude our proof:
combining  (\ref{help_est_0}) and (\ref{help_est}) yields
\begin{eqnarray*}
 \textrm{E} \left| \int_{\Delta^{m}([0,1])} dB^{\alpha}\right|^{2}  \leq
\frac{\gamma_{H}^{|\alpha|}}{\gamma_{H}^{m }} \textrm{E} \left| \int_{\Delta^{m}([0,1])} dB^{(1, \ldots, 1)}\right|^{2}.
\end{eqnarray*}
But clearly $$\int_{\Delta^{m}([0,1])} dB^{(1, \ldots, 1)} = \frac{1}{m!}(B_{1})^{m}$$ and thus we have
$$\textrm{E} \left| \int_{\Delta^{m}([0,1])} dB^{(1, \ldots, 1)}\right|^{2} = \frac{(2m)!}{2^{m}(m!)^{3}}.  $$
Since 
$$ \frac{(2m)!}{2^{m}(m!)^{2}} \leq 2^{m} $$ the assertion (\ref{step1})
follows.
\hfill $\square$

\bigskip

Putting together Propositions \ref{mal_deriv} and \ref{est_int},  we also get the following estimate for the second moment of the Malliavin derivative of an iterated integral.

\bigskip

\begin{proposition}\label{mal_deriv_est}
Let $m \in \N$ and $\alpha \in \mathcal{A}_{m}$. There exists a constant $K_{2}>0$, depending only on $H$ and $T$, such that we have
\begin{eqnarray}{\label{moment:deriv}}
   \left( {\rm E} \left| D_{u}^{i} \int_{\Delta^{m}([s,t])} d B^{\alpha} \right|^{2} \right)^{1/2}  \leq  \left| \mathfrak{J}_{\alpha,i} \right| \frac{K_{2}^{m-1}} {\sqrt{(m-1)!} }  |t-s|^{(|\alpha|-1)H +m - |\alpha|}
 \end{eqnarray}  for $i=1, \ldots, d$ and all $0 \leq s \leq t \leq T$.
\end{proposition}
\noindent{\it Proof.}
Thanks to Proposition  \ref{mal_deriv}  we have that 
\begin{eqnarray*}
&&  D_{u}^{i} \int_{\Delta^{m}([s,t])} d B^{\alpha}    = \sum_{j \in \mathfrak{J}_{\alpha,i}}\int_{   s \leq  t_{1} \leq  \ldots \leq t_{j-1}\leq u \leq t_{j} \leq \ldots \leq t_{m-1}\leq t } \, d B^{\alpha_{-j}}(t_{1}, \ldots, t_{m-1}).  \nonumber
 \end{eqnarray*} Thus it follows
\begin{multline}\label{bnd:deriv-intg2}
\left( {\rm E} \left| D_{u}^{i} \int_{\Delta^{m}([s,t])} d B^{\alpha}(t_{1}, \ldots, t_{m}) \right|^{2} \right)^{1/2}   \\
\leq  \sum_{j \in \mathfrak{J}_{\alpha,i}} \left( {\rm E} \left|\int_{   s \leq  t_{1} \leq  \ldots \leq t_{j-1}\leq u \leq t_{j} \leq \ldots \leq t_{m-1}\leq t } \, d B^{\alpha_{-j}}(t_{1}, \ldots, t_{m-1}) \right|^{2} \right)^{1/2}  
\end{multline}
Furthermore, it is easily checked that
\begin{multline*}
 \int_{   s \leq  t_{1} \leq  \ldots \leq t_{j-1}\leq u \leq t_{j} \leq \ldots \leq t_{m-1}\leq t } \, d B^{\alpha_{-j}}(t_{1}, \ldots, t_{m-1})\\  
=
 \int_{ \Delta^{l(\alpha^{j_{1}})}([s,u])} dB^{\alpha^{j_{1}}} \times \int_{ \Delta^{l(\alpha^{j_{2}})} ([u,t])} dB^{\alpha^{j_{2}}},
  \end{multline*}
with $\alpha=( \alpha^{j_{1}},i,\alpha^{j_{2}})$. Since an iterated integral belongs to a finite chaos with respect to $B$, all its $L^p$ norms are equivalent. See, e.g., Theorem 1.4.1 in \cite{nual}. Thus, we obtain from  Proposition \ref{est_int} and H\"older's inequality that  
\begin{multline}\label{bnd:split}
  \left( {\rm E}\left| \int_{   s \leq  t_{1} \leq  \ldots \leq t_{j-1}\leq u \leq t_{j} \leq \ldots \leq t_{m-1}\leq t } \, d B^{\alpha_{-j}}(t_{1}, \ldots, t_{m-1}) \right|^{2} \right)^{1/2}  \\ 
c_{2,4} |t-s|^{|\alpha_{-j}|H +m-1-|\alpha_{-j}|}  \frac{K_{1}^{m-1}}{\sqrt{l(\alpha^{j_{1}})!}\sqrt{l(\alpha^{j_{2}})!} },
  \end{multline} with a constant $c_{2,4}>0$.
Moreover, it is readily seen that
$$
\frac{1} {\sqrt{l(\alpha^{j_{1}})!}\sqrt{l(\alpha^{j_{2}})!}} \leq  
\frac{1}{[m/2]!},
$$
and according to the fact that $(2k)!/(k!)^2\le 2^{2k}$, we end up with
$$
\frac{1} {\sqrt{l(\alpha^{j_{1}})!}\sqrt{l(\alpha^{j_{2}})!}} \leq  
\frac{2^{(m-1)/2}}{\sqrt{(m-1)!}}.
$$
Plugging this inequality into (\ref{bnd:split}) and (\ref{bnd:deriv-intg2}),
we obtain
\begin{eqnarray*}
&& \left( {\rm E}\left| D_u^i\int_{   s \leq  t_{1} \leq  \ldots \leq t_{j-1}\leq u \leq t_{j} \leq \ldots \leq t_{m-1}\leq t } \, d B^{\alpha_{-j}}(t_{1}, \ldots, t_{m-1}) \right|^{2} \right)^{1/2}  \\ &&  \qquad \quad \qquad \qquad \leq c_{2,4} \frac{(\sqrt{2}K_{1})^{m-1}} {\sqrt{(m-1)!} }  \left| \mathfrak{J}_{\alpha,i} \right| |t-s|^{(|\alpha_{-j}|H +m -1 - |\alpha_{-j}|)},
\end{eqnarray*}
and since
$$ |\alpha_{-j}|H +m -1 - |\alpha_{-j}| =  (|\alpha|-1)H +m - |\alpha|,$$ 
our claim (\ref{moment:deriv}) follows.
\begin{flushright} $\square$ \end{flushright}

\subsection{Study of the remainder term for $H>1/2$}

To avoid notational confusion we will write in the following $X_{t}$, $t \in [0,T]$, instead of 
$X_{t}^{a}$, $t\in [0,T]$, for the solution of the SDE with $X_{0}=a$. Moreover, recall that 
$X^{i}_{t}$, $t \in [0,T]$, denotes the $i$-th component of $X$.
Recall also that the differential operators
$\mathcal{D}^0$ and
$\mathcal{D}^j$  are defined as
\begin{equation}
    \mathcal{D}^0 = \sum_{k=1}^n b^k \, \frac{\partial}{\partial x^k} \qquad \qquad
    \text{ and } \qquad \qquad \mathcal{D}^j = \sum_{k=1}^n \sigma^{k,j} \,
    \frac{\partial}{\partial x^k}
\end{equation}
for $j=1, \ldots, d$ and that we have set
$\mathcal{D}^{\alpha} = \mathcal{D}^{\alpha_1} \ldots
\mathcal{D}^{\alpha_k}$ for a multi-index $\alpha \in \mathcal{A}_{k}$. 

With the help of the auxiliary results contained in the previous section,
we are now able to bound ${\rm E}\mathcal{R}_{m}(0,t)$ in the
following way when $H>1/2$:
\begin{theorem}{\label{remainder}} Let $m \in \N$, $H>1/2$ and assume that  assumption (A) holds. Then there exists a constant $K_{3}>0$, depending only on $H$, $T$, $d$ and $n$, such that
 $$  \left| {\rm E}\mathcal{R}_{m}(0,t) \right| \leq  ( \mathcal{U}_{m+1} + \widetilde{\mathcal{U}}_{m+1}\mathcal{Y}) \frac{ K_{3}^{m} t^{H(m+1)}}{ \sqrt{ m!}} $$  for all $t \in [0,T]$, 
where 
\begin{eqnarray*}
 && \mathcal{U}_{m} = \sup_{\alpha \in \mathcal{A}_{m}}  \sup_{0 \leq t \leq T} \left({\rm E} \left|\mathcal{D}^{\alpha}f(X_t) \right|^{2} \right)^{1/2}, \\
&&\widetilde{\mathcal{U}}_{m} = \max_{i=1, \ldots, n}\sup_{\alpha \in \mathcal{A}_{m}}\sup_{0 \leq t \leq T} \left( {\rm E} \left| \frac{\partial }{ \partial x_{i}} \mathcal{D}^{\alpha} f(X_{t}) \right|^{2}    \right)^{1/2}
\end{eqnarray*}
and
$$ \mathcal{Y} =  \max_{i=1, \ldots, n} \max_{j= 1, \ldots, d} \sup_{0 \leq u \leq s \leq T} \left( {\rm E} \left|D^{j}_{u}X_{s}^{i} \right|^{4} \right)^{1/4}. $$ 

\end{theorem}
\bigskip
Notice then that the second part of Theorem
\ref{thmptfa} is an immediate consequence of the above estimate.

Before we can prove Theorem {\ref{remainder}}, we will need the following proposition, which is a  
straightforward consequence of Proposition 5.2.3 in \cite{nual}, Proposition \ref{prop_X} and the 
properties of iterated integrals of fractional Brownian motion.

\bigskip

\begin{proposition}  \label{exp_int} Let  $m \in \N$, $\alpha \in \mathcal{A}_{m}$,
$g \in C^{2}_{b}(\R^{n};\R)$, and set 
$J_{\al}(s,t)=\int_{\Delta^{m}([s,t])} \, dB^{\alpha}$.  Then it  holds
\begin{eqnarray*}
 && {\rm E} \left( \int_{0}^{t} g(X_{s})  J_{\al}(s,t)  \, dB^{j}_{s} \right) \\ && \qquad = \gamma_{H} {\rm E} \left(  \int_{0}^{t} \int_{0}^{s} \sum_{i=1}^{n}  g_{x_{i}}(X_{s}) \,  J_{\al}(s,t) \, D^{j}_{u}X_{s}^{i} \,  |s-u|^{2H-2} \, du \, ds \right) \\ && \qquad \quad   + \gamma_{H} {\rm E} \left( \int_{0}^{t} \int_{s}^{t}   g(X_{s}) \,  D^{j}_{u}    J_{\al}(s,t)  \, |s-u|^{2H-2} \, du \, ds \right) \end{eqnarray*}
for $t \in [0,T]$ and $j=1, \ldots, d$.
\end{proposition}

\bigskip

We are now ready to prove the main result of this section.

\bigskip

\noindent
{\it Proof of Theorem \ref{remainder}:}
Note that by the proof of Theorem \ref{St-tree-expansion-exact-sol:Wm} we have
$$ \mathcal{R}_{m}(0,t)= \sum_{\alpha \in \mathcal{A}_{m+1}}  \int_{0}^t \int_{0}^{t_{m+1}} 
\ldots \int_{0}^{t_2}
    \mathcal{D}^{\alpha} f(X_{t_1}) \, dB_{t_1}^{\alpha_1} \, dB_{t_2}^{\alpha_2} \ldots 
dB_{t_{m+1}}^{\alpha_{m+1}}.$$
(a) We first consider a single integrand.
By interchanging the order of  integration, which is possible since all integrals are pathwise defined,  we have
\begin{eqnarray*} && \int_{0}^{t} \int_{0}^{t_{m+1}}  \cdots  \int_{0}^{t_{2}}   
\mathcal{D}^{\alpha} f(X_{t_{1}}) \, dB^{\alpha_{1}}_{t_{1}} \, dB^{\alpha_{2}}_{t_{2}} 
 \cdots \, dB^{\alpha_{m+1}}_{t_{m+1}} \\ && \qquad =  
\int_{0}^{t} \int_{t_{1}}^{t} \int_{t_{1}}^{t_{m+1}} \cdots \int_{t_{1}}^{t_{3}}    
\, dB^{\alpha_{2}}_{t_{2}}  \cdots \, dB^{\alpha_{m}}_{t_{m}}\, dB^{\alpha_{m+1}}_{t_{m+1}} 
\, \mathcal{D}^{\alpha} f(X_{t_{1}}) \, dB^{\alpha_{1}}_{t_{1}} \\ && \qquad = 
\int_{0}^{t} \int_{\Delta^{m}([s,t])}   \, dB^{ \alpha_{-1}} \, \mathcal{D}^{\alpha} f(X_{s}) \, 
dB^{\alpha_{1}}_{s} .
\end{eqnarray*}
Recall that
$$ \mathcal{U}_{m} = \sup_{\alpha \in \mathcal{A}_{m}}\sup_{0 \leq t \leq T} \left( {\rm E} |\mathcal{D}^{\alpha} f(X_{t})|^{2}    \right)^{1/2}$$
and
$$ \widetilde{\mathcal{U}}_{m} = \sup_{i=1, \ldots, n}\sup_{\alpha \in \mathcal{A}_{m}}\sup_{0 \leq t \leq T} \left( {\rm E} \left| \frac{\partial }{ \partial x_{i}} \mathcal{D}^{\alpha} f(X_{t}) \right|^{2}    \right)^{1/2}.$$
 If $\alpha_{1}=0$, we  clearly have
\begin{eqnarray} \label{bound_alpha_0}\nonumber
&&   \left| {\rm E} \int_{\Delta^{m+1}([0,t])}   \mathcal{D}^{\alpha} f(X_{t_{1}}) \, dB^{\alpha} (t_{1}, \ldots, t_{m+1})  \right|
\\ && \qquad \qquad \leq  \mathcal{U}_{m+1}  \int_{0}^{t}   \left( {\rm E} \left| \int_{\Delta^{m}([s,t])} \, dB^{\alpha_{-1}}  \right|^{2}\right)^{1/2}  \, ds  \nonumber
\\ && \qquad \qquad \leq  \mathcal{U}_{m+1}  K_{1}^{m} \frac{t^{Hm+1}}{\sqrt{m!}}. 
\end{eqnarray}
 If $\alpha_{1}\neq 0$ we have, according to  Proposition  \ref{exp_int}:
\begin{eqnarray*}
&& \left| \textrm{E} \int_{\Delta^{m+1}([0,t])}   \mathcal{D}^{\alpha} f(X_{t_{1}}) \, dB^{\alpha} (t_{1}, \ldots, t_{m+1})  \right|
 \\ && \qquad \leq  \gamma_{H} \left| \int_{0}^{t} \int_{0}^{s}  \,   \sum_{i=1}^{n}\, {\rm E}  \frac{\partial }{\partial x_{i}} \mathcal{D}^{\alpha} f (X_{s})   \,  \int_{\Delta^{m}([s,t])} \, dB^{\alpha_{-1}} \, D^{\alpha_{1}}_{u}X_{s}^{i}  \,   |s-u|^{2H-2} \, du \, ds \right|\\ && \qquad \quad   +  \gamma_{H}  \left| \int_{0}^{t} \int_{s}^{t}  {\rm E}  \, \mathcal{D}^{\alpha}f(X_{s}) \,\,   D^{\alpha_{1}}_{u}    \int_{\Delta^{m}([s,t])} \, dB^{\alpha_{-1}} \, |s-u|^{2H-2} \, du \, ds  \right|.\end{eqnarray*}
Thus it follows
\begin{eqnarray*}
&& \left| { \rm E} \int_{\Delta^{m+1}([0,t])}   \mathcal{D}^{\alpha} f(X_{t_{1}}) \, dB^{\alpha} (t_{1}, \ldots, t_{m+1})  \right|
 \\ && \,\, \leq  \sum_{i=1}^{n} \widetilde{\mathcal{U}}_{m+1} \gamma_{H} \int_{0}^{t} \int_{0}^{s} \left( {\rm E} \left |D^{\alpha_{1}}_{u}X_{s}^{i} \right|^{4} \right)^{1/4} \left(  {\rm \E } \left| \int_{\Delta^{m}([s,t])} \, dB^{\alpha_{-1}} \right|^{4} \right)^{1/4}  |s-u|^{2H-2} \, du \, ds \\ && \qquad \quad   +  \mathcal{U}_{m+1}\gamma_{H}\int_{0}^{t} \int_{s}^{t}   \left( {\rm \E } \left|  D^{\alpha_{1}}_{u}   \int_{\Delta^{m}([s,t])} \, dB^{\alpha_{-1}} \right|^{2} \right)^{1/2}  |s-u|^{2H-2} \, du \, ds  .\end{eqnarray*}
So recalling that we have set
$$ \mathcal{Y} =  \max_{i=1, \ldots, n} \max_{j= 1, \ldots, d} \sup_{0 \leq u \leq s \leq T} \left( {\rm E} \left |D^{j}_{u}X_{s}^{i} \right|^{4} \right)^{1/4}, $$
we get
\begin{eqnarray*}
&& \left| { \rm E} \int_{\Delta^{m+1}([0,t])}   \mathcal{D}^{\alpha} f(X_{t_{1}}) \, dB^{\alpha} (t_{1}, \ldots, t_{m+1})  \right|
 \\ && \,\, \qquad   \leq  n \widetilde{\mathcal{U}}_{m+1} \mathcal{Y} \gamma_{H} \int_{0}^{t} \int_{0}^{s}  \left(  {\rm \E } \left| \int_{\Delta^{m}([s,t])} \, dB^{\alpha_{-1}} \right|^{4} \right)^{1/4}  |s-u|^{2H-2} \, du \, ds \\ && \qquad \quad   \qquad  +  \mathcal{U}_{m+1}\gamma_{H}\int_{0}^{t} \int_{s}^{t}   \left( {\rm \E } \left|  D^{\alpha_{1}}_{u}   \int_{\Delta^{m}([s,t])} \, dB^{\alpha_{-1}} \right|^{2} \right)^{1/2}  |s-u|^{2H-2} \, du \, ds  .\end{eqnarray*}
Furthermore, invoking again the equivalence of $L^p$ norms
for iterated integral and  Proposition \ref{est_int}, we obtain
\begin{eqnarray*}
 && \gamma_{H} \int_{0}^{t} \int_{0}^{s}  \left(  {\rm \E } \left| \int_{\Delta^{m}([s,t])} \, dB^{\alpha_{-1}} \right|^{4} \right)^{1/4}  |s-u|^{2H-2} \, du \, ds  \\&& \qquad \qquad \leq c_{2,4} \frac{K_{1}^{m}}{\sqrt{m!}}  \gamma_{H} \int_{0}^{t} \int_{0}^{s} |t-s|^{Hm} \, |s-u|^{2H-2} \, du \, ds  \\ && \qquad \qquad \leq c_{2,4} K_{1}^{m} \frac{t^{H(m+2)}}{\sqrt{m!}}.
\end{eqnarray*}
By Proposition \ref{mal_deriv_est}, we get similarly 
\begin{eqnarray*}
&& \gamma_{H}\int_{0}^{t} \int_{s}^{t}   \left( {\rm \E } \left|  D^{\alpha_{1}}_{u}   \int_{\Delta^{m}([s,t])} \, dB^{\alpha_{-1}} \right|^{2} \right)^{1/2}  |s-u|^{2H-2} \, du \, ds  \\ && \qquad \qquad \leq c_{2,4} K_{2}^{m-1}  \frac{t^{H(m+1)}}{\sqrt{(m-1)!}}. \end{eqnarray*}
Thus, we have shown for $\alpha_{1} \neq 0$ the estimate
\begin{eqnarray} \label{bound_alpha_neq0}
&& \left| { \rm E} \int_{\Delta^{m+1}([0,t])}   \mathcal{D}^{\alpha} f(X_{t_{1}}) \, dB^{\alpha} (t_{1}, \ldots, t_{m+1})  \right|  
 \\ && \,\,\qquad \qquad \qquad   \leq  n \widetilde{\mathcal{U}}_{m+1} \mathcal{Y} c_{2,4} K_{1}^{m} \frac{t^{H(m+2)}}{\sqrt{m!}}  + \mathcal{U}_{m+1}c_{2,4}  K_{2}^{m-1} \frac{t^{H(m+1)}}{\sqrt{(m-1)!}}. \nonumber
\end{eqnarray}
(b)  Now we consider the complete remainder term. We have
\begin{eqnarray*}
&& |{\rm E}\mathcal{R}_{m}(0,t)| \leq   \sum_{\alpha \in \mathcal{A}_{m+1}, \alpha_1 =0 }    \left| \textrm{E} \int_{\Delta^{m+1}([0,t])}   \mathcal{D}^{\alpha} f(X_{t_{1}}) \, dB^{\alpha} (t_{1}, \ldots, t_{m+1})  \right|
\\ \nonumber  && \qquad \qquad \qquad \qquad   +
\sum_{\alpha \in \mathcal{A}_{m+1}, \alpha_{1} \neq 0 } \left| \textrm{E} \int_{\Delta^{m+1}([0,t])}   \mathcal{D}^{\alpha} f(X_{t_{1}}) \, dB^{\alpha} (t_{1}, \ldots, t_{m+1})  \right|
\end{eqnarray*}  Since $|\mathcal{A}_{m}| = (d+1)^{m}$,
it follows by (\ref{bound_alpha_0}) and (\ref{bound_alpha_neq0}), that there exists a constant $K_{3}>0$ depending only on $H$, $T$, $n$ and $d$ such that  
 \begin{eqnarray*}
&& |{\rm E}\mathcal{R}_{m}(0,t)| \leq   \left( \mathcal{U}_{m+1} + \mathcal{Y} \widetilde{\mathcal{U}}_{m+1}\right) K_{3}^{m} \frac{t^{H(m+1)}}{\sqrt{m!}},
\end{eqnarray*} 
which completes the proof.
\begin{flushright} $\square$ \end{flushright}

\bigskip

\noindent
{\bf Acknowledgments}. Part of this work was done during the stays of authors at our respective institutions
and we would like to thank them for their support and hospitality.

\end{document}